# A Simple Approach to the Tiling Problem Using Recursive Sequence

*A Project Submitted for the Singapore Mathematics Olympiad Project Festival 2019*

**Done by:**

Le Viet Hung

Jin Qingyang

Huang Keyi

Tan Yiming

ANDERSON SECONDARY SCHOOL

**Project Consultant:** Dr Chua Seng Kiat

**Project Supervisor:** Mr Mohamed Alfin H R

# Table of contents

## *Abstract*

The tiling problem has been a famous problem that has appeared in many Mathematics problems. Many of its solutions are rooted in high level Mathematics. Thus we hope to tackle this problem using more elementary Mathematics concepts. In this report, we start with the simplest cases, with the smaller numbers: the number of ways to tile a 2 x *n*, 3 x *n*, 4 x *n* rectangular board using 2 x 1 domino tiles, where the number of rows are fixed and we present a recursive formula based on m and the earlier terms. This allows us to deduce the non-recursive formula for each case that is only dependent on *m*. For each case, we also expand and generalize the problem, not just for 2,3,4 but for any positive integer *k*, for certain types of configurations of the board. We also focus on 1 of the famous variations of the tiling problem: tatami tiling, and present a solution for simple cases: 2 x *n*, 3 x *n*, 4 x *n*. In the end, we have managed to find a simpler solution for 3 different configurations of the board, with some we even deduced the non-recursive formula. We have also solved simple cases of the tatami tiling problem, with hope to tackle the general case in the future. We realized that our method only works on a case-by-case basis, with little success in solving the general case. This approach is also applicable in many other counting problems which we wish to pursue for further research.

# I.   Introduction

Tiling of rectangular boards is a famous concept that appears in many mathematics problems, especially in the field of Combinatorics. On the simplest level, the question appears in many math problems: "Is it possible to tile a specific-shape board X with Y tiles of size …?" (for example, the mutilated chessboard problem). However, in cases when we have established that it is possible, we ask ourselves the next, more difficult question: "How many ways are there?".

Of course this question is no easy feat for weird board shapes and weird-shaped tiles. Thus here we just focus on the simplest case: tiling a rectangular board with 1 x *m* rectangular tiles. For example: how many ways are there to tile a *m* x *n* rectangular board with 2 x1 dominoes?

This general problem has been solved for all cases of *m* x *n*. In 1961, Temperley & Fisher and Kasteleyn both independently calculated a beautiful formula:

$$\prod_{j=1}^{m}\prod_{k=1}^{n}\sqrt[4]{4\cos^2\frac{\pi j}{m+1}+4\cos^2\frac{\pi k}{n+1}}$$

This formula is applicable to all cases of *m* and *n*. However, the formula involves certain trigonometric terms and its proof involves very high level algebra that may not be understandable to most students.

Thus we wanted to work on this project to find a simpler approach to tackling the problem that only involves lower, pre-university level mathematics knowledge. We were interested in number patterns and found a way to apply these concepts in the tiling question. And the most intuitive method is through recursive sequence. This project also helped us to widen our Mathematics knowledge and skills, especially in the fields of Combinatorics, Algebra and simple Calculus.

Another area of our interest is the tatami tiling problem. It is the same as the regular tiling problem, with a twist: each possible tatami tiling is subject to the constraint that no four corners meet at a point. Again, this problem has been beautifully solved, however its solution is based on generating functions. Thus we wish to apply our approach to produce a simpler solution.

## II. Domino tiling of 2 x *m* rectangular board and its variations
### 1. Domino tiling of 2 x *m* rectangular board

In this report, generally our approach for all problems is to tile sequentially from the left to the right of the board, examining each case extensively to assure that all cases have been covered and no two cases overlap.

We will begin with the simplest case: the number of ways to tile an 2 x *m* rectangular board using 2x1 tiles, a(*m*).

**Solution**

For *m* = 1, it is obvious that there is only 1 way to tile a 2 x 1 rectangular board using one 2 x 1 tile. Thus **a(1)=1**.

For *m* = 2, we have a 2 x 2 board. Clearly there are only 2 ways to tile the board (either horizontally or vertically). Thus **a(2)=2**.

For *m*≥3, consider 2 cases:

a) When the tile covering the bottom left corner of the board is vertical, we are left with a 2 x (*m*-1) board to tile. **There are a(*m*-1) ways to tile this board.**

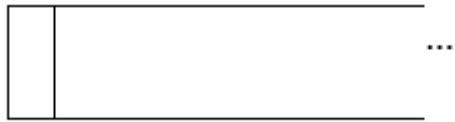

b) When the tile covering the bottom left corner of the board is horizontal, then the tile covering the top left corner is also horizontal. We are left with a 2 x (*m*-2) board to tile. **There are a(*m*-2) ways to tile this board.**

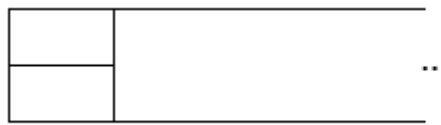

To summarize, the number of ways to tile a 2 x *m* rectangular board using 2 x 1 tiles,

**a(*m*)=a(*m*-1)+a(*m*-2)** (for *m*≥3)

with a(1)=1 and a(2)=2

We will now attempt to find the non-recurrence formula for a (*m*), such that it is only dependent on *m* and not on the earlier terms of the sequence.

The characteristic equation is

$$x^2 = x + 1$$
$$\Leftrightarrow x^2 - x - 1 = 0$$
$$\Leftrightarrow x = \frac{1 \pm \sqrt{5}}{2}$$

The non-recurrence formula of a(m) has the form

$$a(m) = ax_1^m + bx_2^m$$

Sub $m = 1$: $a(1) = ax_1 + bx_2$

$$\Rightarrow 2 = a(1 + \sqrt{5}) + b(1 - \sqrt{5}) \quad (1)$$

Sub $m = 2$: $a(2) = ax_1^2 + bx_2^2$

$$\Rightarrow 2 = a(\frac{1+\sqrt{5}}{2})^2 + b(\frac{1-\sqrt{5}}{2})^2$$

$$\Rightarrow 8 = a(6 + 2\sqrt{5}) + b(6 - 2\sqrt{5}) \quad (2)$$

From (1) and (2):

$$a = \frac{5+\sqrt{5}}{10} \text{ and } b = \frac{5-\sqrt{5}}{10}$$

Thus the non-recursive formula for a(m) is

$$a(m) = (\frac{5+\sqrt{5}}{10})(\frac{1+\sqrt{5}}{2})^m + (\frac{5-\sqrt{5}}{10})(\frac{1-\sqrt{5}}{2})^m$$

It is interesting to note that the recursive formula for a(*m*) is also the formula for the familiar Fibonacci sequence, and all the terms in a(*m*) are also terms of the Fibonacci sequence.

## 2. Tromino tiling of 3 x *m* rectangular board

The domino tiling problem of the 2 x *m* board was simple enough, and it yields a rather beautiful solution that is closely related to the Fibonacci sequence. The next question in our minds would be if it is possible to use this approach for the tiling by trominoes of a 3 x *m* board?

It turns out that we can. Indeed,

**Solution**

For *m* = 1, it is obvious that there is only 1 way to tile a 3 x 1 rectangular board using one 3 x 1 tile. Thus **a(1)=1**.

For m = 2, we have a 3 x 2 board. Clearly there is only 1 way to tile the board, using 2 vertical 3 x 1 trominoes. Thus **a(2)=1**.

For m = 3, we have a 3 x 3 board. There are 2 ways to tile the board (either by using 3 horizontal or 3 vertical trominoes). Thus **a(3)=2**.

For *m*≥4, consider 2 cases involving the tiling of the leftmost column:

a) When the leftmost column is covered by 1 vertical 3x1 tromino, we are left with a 3 x (*m*-1) board to tile. **There are a(*m*-1) ways to tile this board.**

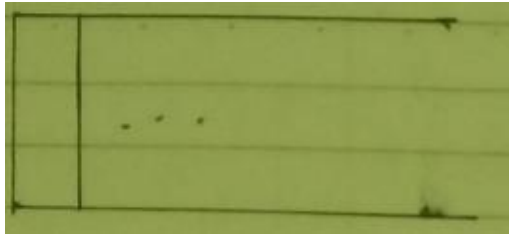

b) When the leftmost column is covered by 3 horizontal 3x1 trominoes, we are left with a 3 x (*m*-3) board to tile. **There are a(*m*-3) ways to tile this board.**

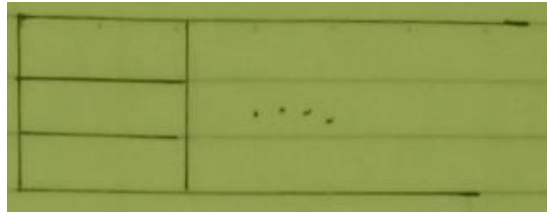

To summarize, the number of ways to tile a 3 x *m* rectangular board using 3 x 1 trominoes,

**a(m)=a(m-1)+a(m-3)** (for *m*≥4)

with a(1)=a(2)=1 and a(3)=2

And here it is possible to make a nice observation that the formula has the exact same form as the one above for the 2 x *m* case. This leads us to the question: Does it work the same way for 4 x *m*, 5 x *m*? Does the solution depends on the number 2,3,4,5, or we can just substitute it with any random positive integer *k*?

## 3. Tiling of *m* x *n* rectangular board by *m* x 1 tiles

We figured out that this method also works for the general case, however, there are some slight modifications to the trivial cases.

**Solution**

Using a similar approach, let's consider the number of ways to tile a *m* x *n* rectangular board using 1 x *m* tiles, a(*m,n*) (with *m*≥2).

**For *n*=1**, it is obvious that there is only 1 way to tile the *m* x 1 rectangular board. Thus **a(m,1)=1**.

**For 2≤*n*≤*m*-1** (*m*≥3), there is also only 1 way to tile the *m* x *n* rectangular board, since no tiles can be horizontal. Thus **a(*m,n*)=1** for 2≤*n*≤*m*-1 (*m*≥3)

**For *n*=*m***, we have a *m* x *m* board. There are 2 ways to tile the board (either all the tiles are horizontal or vertical). Thus **a(*m,m*)=2**.

**For *n*≥*m*+1**, consider 2 cases:

a) When the tile covering the bottom left corner is vertical, then for the remaining *m* x (*n*-1) board there are a(*m,n*-1) ways to tile.

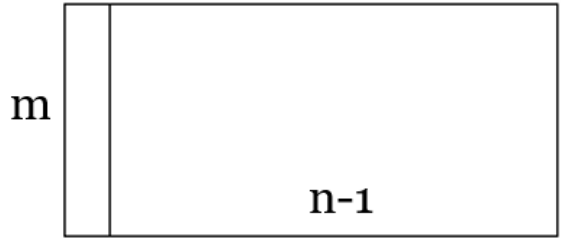

b) When the tile covering the bottom left corner is horizontal, then the other tiles are "forced" into place as follows.

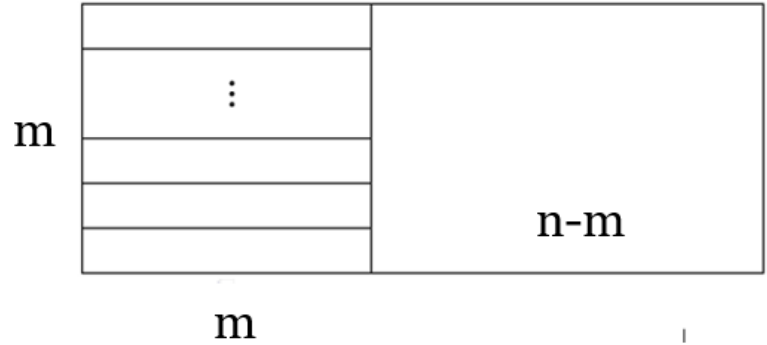

We are left with a *m* x (*n-m*) board, and there are a(*m,n-m*) ways to tile.

In total there are a(*m,n*-1)+a(*m,n-m*) ways to tile. Thus in this case **a(*m,n*)= a(*m,n*-1)+a(*m,n-m*).**

To summarize, the number of ways to tile a *m* x *n* rectangular board using *m* x 1 tiles is

**a(*m*,1)=a(*m*,2)=a(*m*,3)=…=a(*m,m*-1)=1; a(*m,m*)=2**   and

**a(*m,n*)= a(*m,n*-1)+a(*m,n-m*)**   for $n \geq m+1$

## 4. Table of values

Using the formula above, we can construct the table of values for the general case of a(*m,n*), the number of tilings of a *m* x *n* rectangular board using *m* x 1 tiles, for relatively small values of *m* and *n*. In this table, the first column represents values of *m*, while the first row represent values of *n*.

|    | 1 | 2 | 3 | 4 | 5 | 6  | 7  | 8  | 9  | 10 | 11  | 12  | 13  | 14  | 15  |
|----|---|---|---|---|---|----|----|----|----|----|-----|-----|-----|-----|-----|
| 2  | 1 | 2 | 3 | 5 | 8 | 13 | 21 | 34 | 55 | 89 | 144 | 233 | 377 | 610 | 987 |
| 3  | 1 | 1 | 2 | 3 | 4 | 6  | 9  | 13 | 19 | 28 | 41  | 60  | 88  | 129 | 189 |
| 4  | 1 | 1 | 1 | 2 | 3 | 4  | 5  | 7  | 10 | 14 | 19  | 26  | 36  | 50  | 69  |
| 5  | 1 | 1 | 1 | 1 | 2 | 3  | 4  | 5  | 6  | 8  | 11  | 15  | 20  | 26  | 34  |
| 6  | 1 | 1 | 1 | 1 | 1 | 2  | 3  | 4  | 5  | 6  | 7   | 9   | 12  | 16  | 21  |
| 7  | 1 | 1 | 1 | 1 | 1 | 1  | 2  | 3  | 4  | 5  | 6   | 7   | 8   | 10  | 13  |
| 8  | 1 | 1 | 1 | 1 | 1 | 1  | 1  | 2  | 3  | 4  | 5   | 6   | 7   | 8   | 9   |
| 9  | 1 | 1 | 1 | 1 | 1 | 1  | 1  | 1  | 2  | 3  | 4   | 5   | 6   | 7   | 8   |
| 10 | 1 | 1 | 1 | 1 | 1 | 1  | 1  | 1  | 1  | 2  | 3   | 4   | 5   | 6   | 7   |

# III. Domino tiling of 3 x *m* rectangular board and its variations
## 1. Domino tiling of 3 x *m* rectangular board

Let's move on from the previous case to a slightly more challenging problem: how many ways are there to tile a 3 x *m* rectangular board with 2 x 1 dominoes?

Obviously this is not as easy as the previous problems as the configuration of the board with relative to the tile has changed. However, by slight modification to our method, we are still able to arrive at the solution for the problem. From here, we will consider the orientation of certain tiles: whether it is horizontal or vertical, and each would yield a separate case. The ideas of "forcing" tiles also emerges.

**Solution**

Let the number of ways to tile a 3 x *m* rectangular board (3 rows, *m* columns) by 2 x 1 dominoes be a(*m*).

It is obvious that for any possible tiling to occur, 3 x *m* must be divisible by 2, hence *m* must be even. Thus **a(*m*)=0 for odd *m***.

We consider the case when *m* is even: let *m*=2*k* (*k* is an integer).

**For *k* = 1**, we have a 3 x 2 board, and there are 3 ways to tile this board with 2 x 1 dominoes:

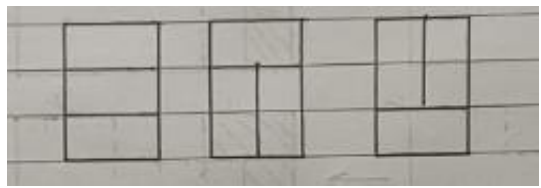

Thus **a(2)=3**.

**For *k*=2**, we have a 3 x 4 board, we consider 3 cases involving the leftmost column:

a) When the leftmost column is covered with 3 horizontal 2x1 dominoes:

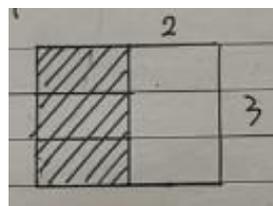

The remaining untiled portion is the same as a 3x2 board. There are a(2) ways to tile.

b) When the leftmost column is covered with 1 horizontal and 1 vertical 2 x 1 domino forming a "⌐" shape:

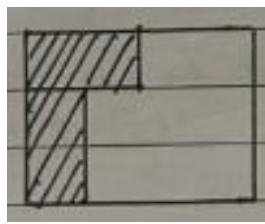

We continue to consider 2 subcases:
- If the third tile is placed as follows, then we are left with a 3 x 2 board, and there are a(2) ways to tile.

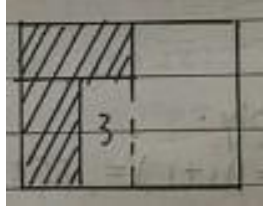

- If the third tile is placed as follows, then the remaining tiles are "forced" into place, and there is only 1 way to tile.

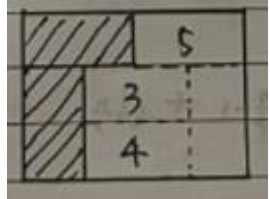

Thus, for case b, in total there are a(2)+1 ways to tile the board.

c) When the leftmost column is covered with 1 horizontal and 1 vertical 2 x 1 domino forming an "L" shape:

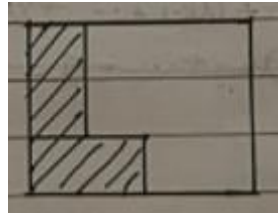

This case is symmetrical to case b, hence there are also a(1)+1 ways to tile the board.

Thus, in total **a(4)=a(2)+2(a(2)+1)=11**.

**For k≥3**, using the same approach, we consider 3 cases involving the tiling of the leftmost column:

a) When the leftmost column is covered with 3 horizontal 2x1 dominoes, then we are left with a 3 x (*m*-2) board to tile, and there are a(*m*-2) ways to tile.
b) When the leftmost column is covered with 1 horizontal and 1 vertical 2x1 domino forming a "Γ" shape:

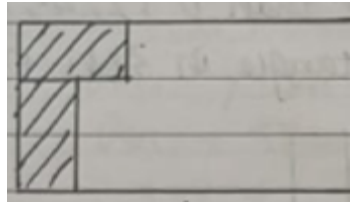

We continue to consider 2 subcases:
- If the third tile is placed as follows, then we are left with a 3 x (*m*-2) board, and there are a(*m*-2) ways to tile.

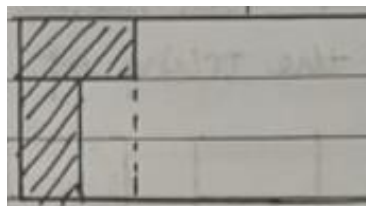

- If the third tile is placed as follows, then the next tiles are "forced":

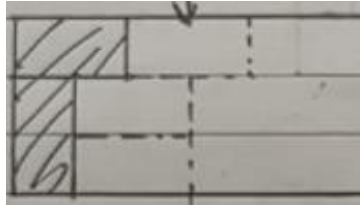

- If the next tile is tiled vertically as follows, then we are left with a 3x(*m*-4) board, and there are a(*m*-4) ways to tile.

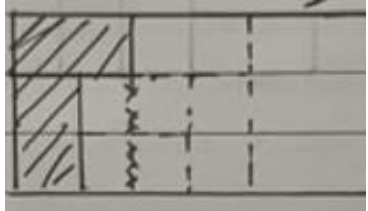

- If the next tile is tiled horizontally as follows, then the next tiles are "forced":

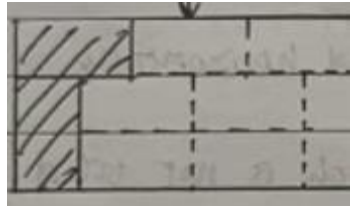

We will repeat this process in a "zig zag" manner until we reach the end of the 3 x *m* board, when the last tile is "forced":

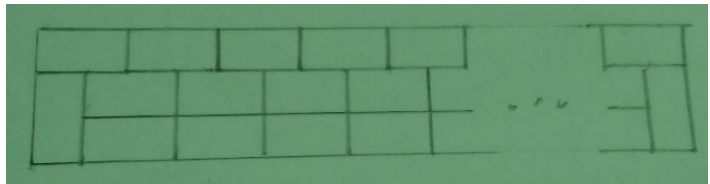

Thus for case **b**, there are **a(*m*-2)+a(*m*-4)+…+a(2)+1** ways to tile the board.

c) When the leftmost column is covered with 1 horizontal and 1 vertical 2 x 1 domino forming an "L" shape:

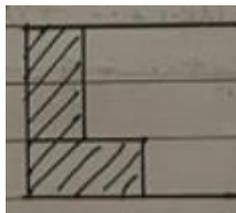

This case is symmetrical to case **b**, thus similarly there are **a(*m*-2)+a(*m*-4)+…+a(2)+1** ways to tile the board.

To summarize, the total number of ways to tile the board is

a(*m*)=a(*m*-2)+2(a(*m*-2)+a(*m*-4)+…+a(2)+1)=**3a(*m*-2)+2a(*m*-4)+2a(*m*-6)+…+2a(2)+2**

Thus a(*m*)= 3a(*m*-2) + 2a(*m*-4) + 2a(*m*-6) + … + 2a(2) + 2

a(*m*-2) = 3a(*m*-4) + 2a(*m*-6) + … + 2a(2)+2

⇨ a(*m*) – a(*m*-2)= 3a(*m*-2) – a(*m*-4)
⇨ **a(*m*) = 4a(*m*-2) - a(*m*-4)** with a(2)=3 and a(4)=11.

The formula also holds true for odd *m*, since a(1) = a(3) = … = a(*m*) = 0.

We will now construct a non-recurrence formula for a(*m*) when *m* is even.

We have a(*2k*) = 4a(*2k*-2) – a(*2k*-4).

Let another sequence b(*k*) be defined as b(*k*) = a(*2k*)=a(*m*), then

b(1)=3, b(2)=11 and **b(*k*) = 4b(*k*-1) – b(*k*-2)**

The characteristic equation is

$$x^2 = 4x - 1$$
$$\Leftrightarrow x^2 - 4x + 1 = 0$$
$$\Leftrightarrow x = 2 \pm \sqrt{3}$$

The non-recurrence formula of b(k) has the form

$$b(k) = bx_1^k + bx_2^k$$

*Sub* $m = 1$: $b(1) = ax_1 + bx_2$

$$\Rightarrow 3 = a(2+\sqrt{3}) + b(2-\sqrt{3}) \qquad (1)$$

*Sub* $m = 2$: $b(2) = ax_1^2 + bx_2^2$

$$\Rightarrow 11 = a(2+\sqrt{3})^2 + b(2-\sqrt{3})^2$$
$$\Rightarrow 11 = a(7+4\sqrt{3}) + b(7-4\sqrt{3}) \qquad (2)$$

*From* (1) *and* (2):

$$a = \frac{3+\sqrt{3}}{6} \text{ and } b = \frac{3-\sqrt{3}}{6}$$

Thus the non-recursive formula for b(k) is

$$b(k) = (\frac{3+\sqrt{3}}{6})(2+\sqrt{3})^k + (\frac{3-\sqrt{3}}{6})(2-\sqrt{3})^k$$

$$\Leftrightarrow a(m) = (\frac{3+\sqrt{3}}{6})(2+\sqrt{3})^{\frac{m}{2}} + (\frac{3-\sqrt{3}}{6})(2-\sqrt{3})^{\frac{m}{2}}$$

when *m* is odd.

## 2. Tromino tiling of 4 x *m* rectangular board

At this stage, it is intuitive to repeat what we did for the earlier case: to replace 2 by 3 and see if the method would still work.

The question now is: how many ways are there to tile a 4 x *m* rectangular board with 3 x1 trominoes?

Obviously substituting 2 by 3 would greatly increase the level of complexity and introduce many more cases compared to the domino tiling of a 3 x *m* rectangular board. Thus here we will introduce and apply a new "technique".

**Solution**

Using a similar approach, let's discover the number of ways to tile a 4 x *m* rectangular board (4 rows, *m* columns) using 3 x 1 trominoes, a(*m*).

It is obvious that for any possible tiling to occur, 4 x *m* must be divisible by 3, hence *m* must be divisible by 3. Thus **a(*m*)=0 for *m* not divisible by 3**.

We consider the case when *m* is divisible by 3: *m*=3*k* (*k* is an integer).

Since there are too many cases to be considered such as in the case for a 3 x *m* board, we will use a slightly different approach to count the number of ways of tiling, while still ensuring that no 2 cases overlap each other:

We start tiling from the left to the right of the board. Draw *k*-1 vertical lines to separate the 4 x *m* board into *m* smaller 4 x 3 boards. For any proper tiling, if the line drawn cuts through the middle of any horizontal tile, we say that the line "intersects" that tile. If there exists a line which doesn't intersect any tile, then to the right of that line remains an untiled portion of a 4 x (3*n*) board, and there are a(3*n*) ways to tile that portion with 3 x 1 trominoes.

**For *k*=1**, we have a 4 x 3 board, and there are 3 ways to tile this board with 3 x 1 trominoes:

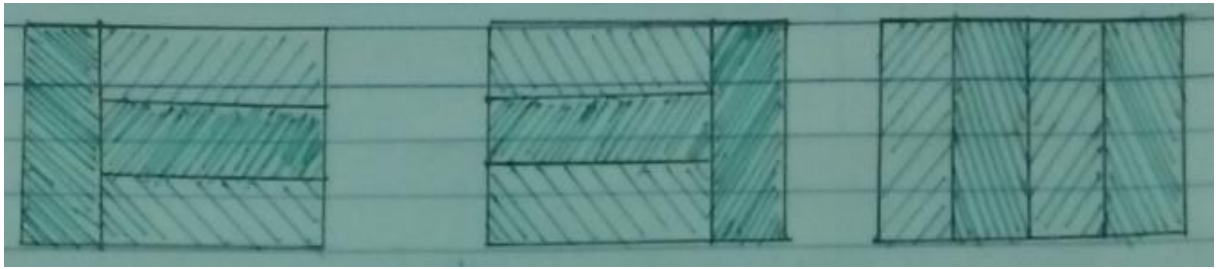

Thus **a(3)=3**.

**For *k*=2**, we have a 4 x 6 board, we consider 3 cases involving the leftmost column:

d) When the leftmost column is covered with 4 horizontal 3x1 trominoes, the remaining untiled portion is the same as a 4x3 board. There are a(3) ways to tile.
e) When the leftmost column is covered with 1 horizontal and 1 vertical 3x1 tromino forming a " ⌐" shape, we continue to consider 3 subcases:
  - If the third and fourth tiles are placed as follows, then we are left with a 4x3 board, and there are a(3) ways to tile.

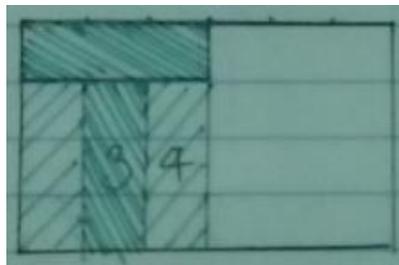

  - If the third and fourth tiles are placed as follows, then the remaining tiles are "forced" into place, and there is only 1 way to tile.

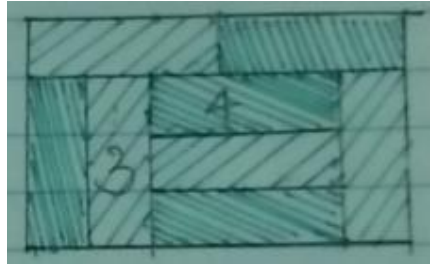

- If the third tile is placed as follows, then the remaining tiles are "forced" into place, and there is only 1 way to tile.

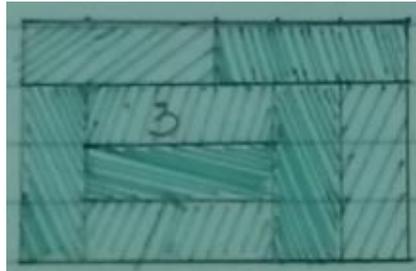

Thus for case b, in total, there are a(3)+2 ways to tile the board.

f) When the leftmost column is covered with 1 horizontal and 1 vertical 3x1 tromino forming an "L" shape, this case is symmetrical to case b, hence similarly there are also a(3)+2 ways to tile the board.

Thus, in total, **a(6)=a(3)+2(a(3)+2)=3a(3)+4=13**.

**For $k=3$**, we have a 4x9 board, using the same approach, we consider 3 cases involving the tiling of the leftmost column:

a) When the leftmost column is covered with 4 horizontal 3x1 trominoes, the remaining untiled portion is the same as a 4x6 board. There are a(6) ways to tile.

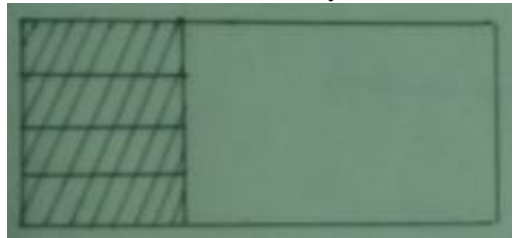

b) When the leftmost column is covered with 1 horizontal and 1 vertical 3x1 tromino forming a "⌐" shape:
As stated above, draw 2 lines separating the 4x9 board into 3 smaller 4x3 boards. Starting from the line on the left, we consider 3 cases:
- The line on the left doesn't intersect any tile, as follows. The remaining untiled portion is a 4x6 board, and there are a(6) ways to tile the board.

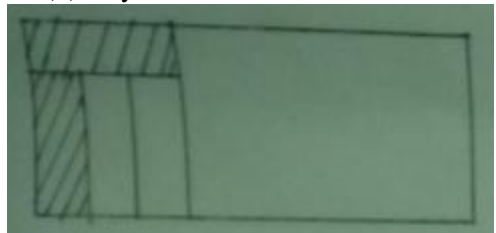

- The line on the left is intersected by a tile while the line on the right doesn't intersect any tile, for example, as follows:

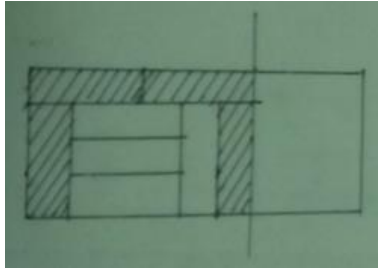

Consider the side of the board to the left of the line on the right:

In the bottom three throws, it is obvious that if any horizontal tile exists, it must exist in groups of 3 forming 3x3 squares. Since the line on the left is intersected by a group of 3 horizontal tiles, the next tile in the top row must be horizontal, since there wouldn't be enough space for a vertical tile.

We can also deduce that the tile on the $6^{th}$ column, bottom row must be vertical (if it is horizontal, because the line on the right is not intersected, so the tile must exist in a 3x3 square stretching from the $4^{th}$ to $6^{th}$ column, which means that the line on the left is not intersected – contradiction).

Thus, the remaining untiled space must consist of 1 3 x 3 square (made up of 3 horizontal tiles) and 1 vertical tile. There are 2 ways to arrange this:

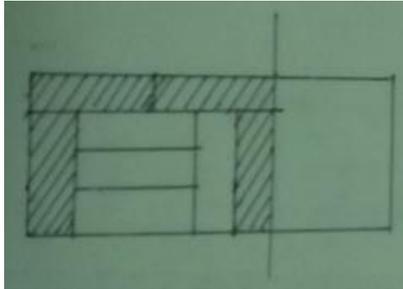 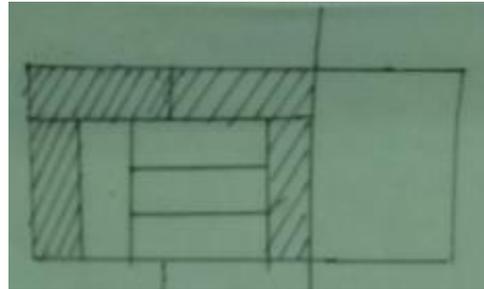

Lastly, on the side of the board to the right of the line, the remaining untiled portion is a 4x3 board, and there are a(3) ways to tile the board.

In conclusion, for this case, there are **2a(3) ways** to tile the board.

- Both lines are intersected:

Using the same reasoning, we can deduce that the top row is filled with 3 horizontal tiles, and the tile on the rightmost column is vertical. Hence the board will be as follows:

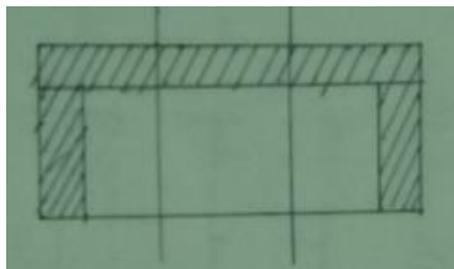

Since the 2 lines are 3 columns apart, no 3x1 horizontal tromino can intersect both lines. Hence each line is intersected by a different square of 3 trominoes. Thus the remaining untiled portion consists of 2 3x3 squares and 1 vertical 3x1 tromino. There are 3 ways arrange 2 squares and 1 vertical tromino (the tromino can be placed in either 1 of the 3 gaps created by the 2 squares). One example is as follows:

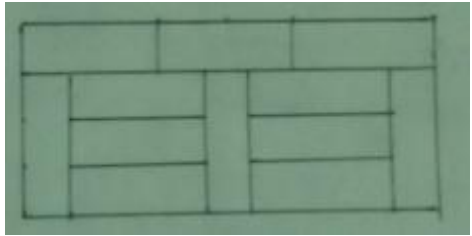

Hence there are 3 ways to tile the board.

In summary, for case **b**, there are **a(6)+2a(3)+3** ways to tile the board.

c) When the leftmost column is covered with 1 horizontal and 1 vertical 3 x 1 tromino forming an "L" shape, this case is symmetrical to case **b**, hence there are also a(6)+2a(3)+3 ways to tile the board.

Thus, **a(9)=a(6)+2(a(6)+2a(3)+3)=3a(6)+4a(3)+6=57**.

For $k \geq 4$, we consider 3 cases involving the tiling of the leftmost column:

a) When the leftmost column is covered with 4 horizontal 3x1 trominoes, the remaining untiled portion is the same as a 4x(*m*-3) board. There are a(*m*-3) ways to tile.
b) When the leftmost column is covered with 1 horizontal and 1 vertical 3x1 tromino forming a "⌐" shape:
   As stated above, draw *k*-1 lines separating the 4 x *m* board into *k* smaller 4 x 3 boards. Number the lines 1, 2, 3, …, *k*-1 starting from the left. Consider the cases:
   - Line 1 doesn't intersect any tile, as follows.

   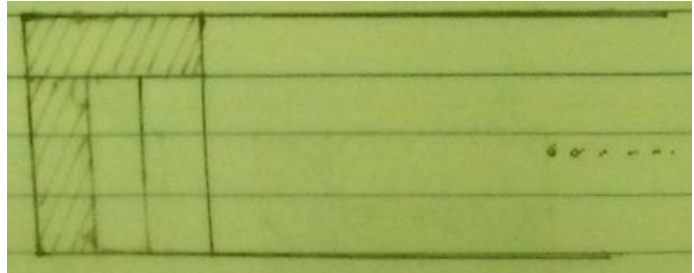

   The remaining untiled portion is a 4x(*m*-3) board, and there are a(*m*-3) ways to tile.
   - Line 1 is intersected while line 2 doesn't intersect any tile. Similar to the above case, there are 2a(*m*-6) ways to tile.
   - We continue this process: suppose line *n* is the first line that doesn't intersect any tile. Then by the same reasoning as above, we can deduce that the top row to the left of line *n* is filled with *n* horizontal 3x1 trominoes, and the tile to the immediate left of line *n* is vertical, as follows.

   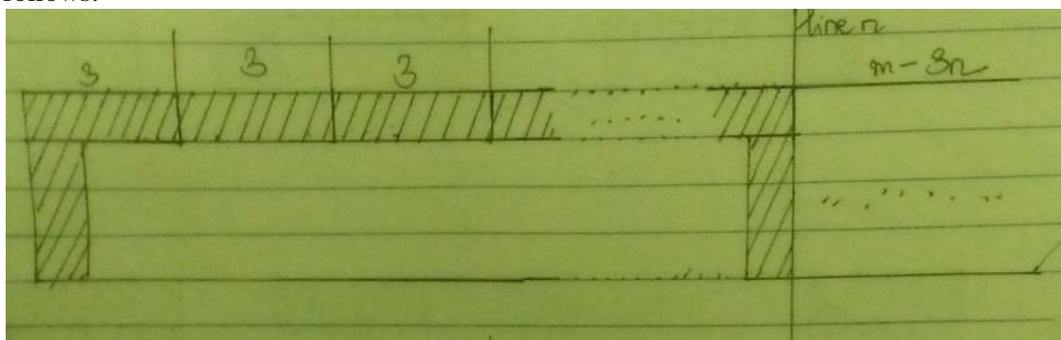

   For the remaining untiled portion to the left of line *n*, since each line from 1 to *n*-1 is intersected by a different 3 x 3 square, hence the whole portion consists of *n*-1 3 x 3 squares and 1 vertical 3 x 1 tromino. There are *n* ways to place the tromino in 1 of the *n* gaps created by the *n*-1 3x3 squares. Hence there are *n* ways to tile this portion.

For the untiled portion to the right of line *n*, it is a 4 x (*m*-3*n*) board, and there are a(*m*-3*n*) ways to tile.

Thus, in total there are **nxa(*m*-3*n*)** ways to tile the board.

- To account for all the cases, then we must consider all the cases of *n* from 1 to *k*-1.
- The last case: when every line is intersected, then the top row is filled with *k* horizontal 3 x1 trominoes, and the rightmost tile is vertical, as follows.

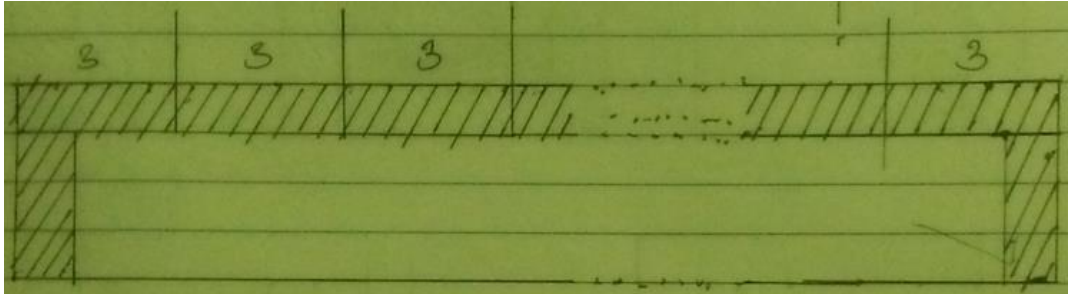

The remaining untiled portion consists of (*k*-1) 3 x 3 squares and 1 vertical 3 x 1 tromino. There are *k* ways to place the tromino in the *k* gaps created by the (*k*-1) squares, thus there are *k* ways to tile the board.

In summary, for case **b**, there are **a(*m*-3)+2a(*m*-6)+3a(*m*-9)+…+(*k*-1)a(3)+*k*** ways to tile the board.

d) When the leftmost column is covered with 1 horizontal and 1 vertical 3 x 1 tromino forming an "L" shape, this case is symmetrical to case **b**, hence there are also **a(*m*-3)+2a(*m*-6)+3a(*m*-9)+…+(*k*-1)a(3)+*k*** ways to tile the board.

Thus, a(*m*)=a(*m*-3)+2(a(*m*-3)+2a(*m*-6)+3a(*m*-9)+…+(*k*-1)a(3)+*k*)

$\quad\quad$ =3a(*m*-3) + 4a(*m*-6) + 6a(*m*-9) + 8a(*m*-12) + … + 2(*k*-1)a(3)+2*k*

$\quad$ a(*m*-3) = 3a(*m*-6) + 4a(*m*-9) + 6a(*m*-12) + … + 2(*k*-2)a(3)+ 2(*k*-1)

$\quad\quad$ a(*m*-6) = 3a(*m*-9) + 4a(*m*-12) + … + 2(*k*-3)a(3) + 2(*k*-2)

⇨ a(*m*)-a(*m*-3)= 3a(*m*-3) + a(*m*-6) + 2a(*m*-9) + 2a(*m*-12) + … + 2a(3) + 2

and $\quad\quad$ a(*m*-3)-a(*m*-6) = 3a(*m*-6) + a(*m*-9) + 2a(*m*-12) + … + 2a(3) + 2

⇨ a(*m*)-a(*m*-3)-[a(*m*-3)-a(*m*-6)]=3a(*m*-3)-2a(*m*-6)+a(*m*-9)

⇨ **a(*m*) = 5a(*m*-3) – 3a(*m*-6) + a(*m*-9)** with a(3)=3, a(6)=13, a(9)=57

This formula also holds true for *m* not divisible by 3, since a(1)=a(2)=a(4)=a(5)=…=a(*m*)=0.

## 3. Tiling of (*k*+1) x *m* rectangular board by *k* x 1 tiles

Obviously substituting 2 by 3 causes a great leap in difficulty level. And here we face the same generalization question again: is it possible to substitute 3 by any random positive integer *k*?

Fortunately, the case for *k* is not very different from the case for 3. However midway we encounters a mini-problem that needs to be resolved, and the result we obtain is much more complex compared to the previous individual cases.

**Solution**

We will now try to expand this problem into the general case: to find the number of tilings of a $(k+1) \times m$ rectangular board ($k+1$ rows, $m$ columns) by $k \times 1$ tiles, $a(k,m)$.

It is obvious that for any possible tiling to occur, $(k+1) \times m$ must be divisible by $k$. Since $GCD(k,k+1)=1$, hence $m$ must be divisible by $k$. Thus **$a(k,m)=0$ for $m$ not divisible by $k$**.

We consider the case when $m$ is divisible by 3: $m=pk$ ($k$ is an integer).

The approach is the same as the above case: to draw $p-1$ vertical lines dividing the board into $p$ $(k+1) \times k$ smaller rectangular boards.

**For $p=1$**, we have a $(k+1) \times k$ rectangular board, and there are 3 ways to tile this board as follows. Thus **$a(k,k)=3$**.

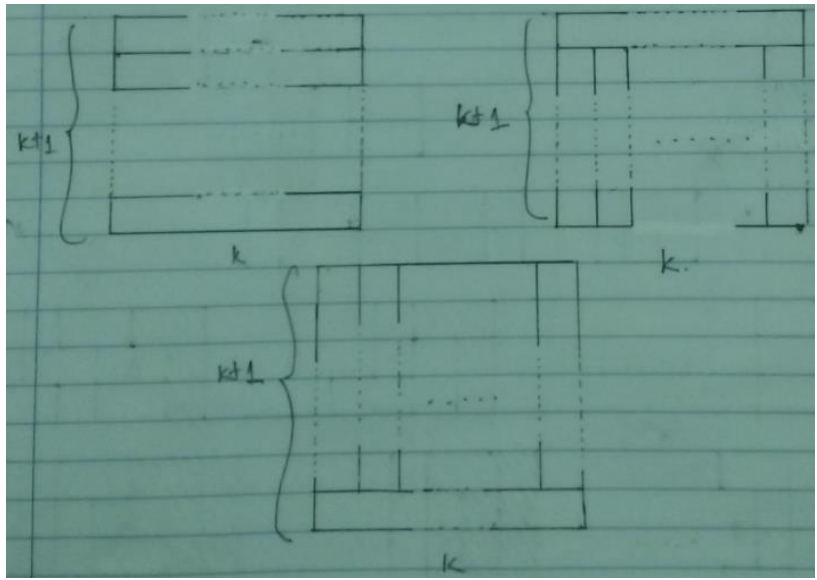

**For $p \geq 2$**, we consider 3 cases involving the tiling of the leftmost column:

c) When the leftmost column is covered with $k+1$ horizontal $k \times 1$ tiles, the remaining untiled portion is the same as a $(k+1) \times (m-k)$ board. There are $a(k,m-k)$ ways to tile.

d) When the leftmost column is covered with 1 horizontal and 1 vertical $k \times 1$ tile forming a "Γ" shape:
As stated above, draw $p-1$ lines separating the $(k+1) \times m$ board into $p$ smaller $(k+1) \times k$ boards. Number the lines 1, 2, 3, …, $p-1$ starting from the left. Consider the cases:
- Line 1 doesn't intersect any tile, as follows.

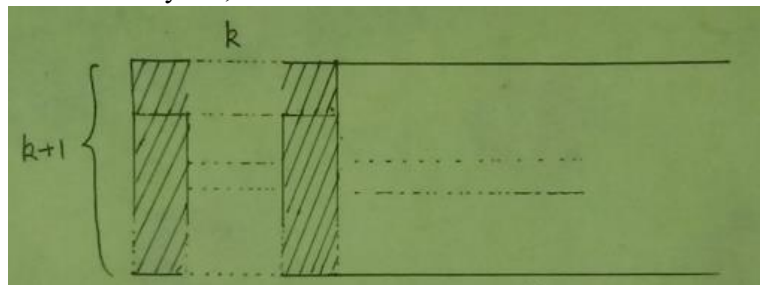

The remaining untiled portion is a $(k+1) \times (m-k)$ board, and there are $a(k,m-k)$ ways to tile.
- We continue this process: suppose line $n$ is the first line that doesn't intersect any tile. Then by the same reasoning as above, we can deduce that the top row to the left of line $n$ is filled with $n$ horizontal $k \times 1$ tiles, and the tile to the immediate left of line $n$ is vertical, as follows.

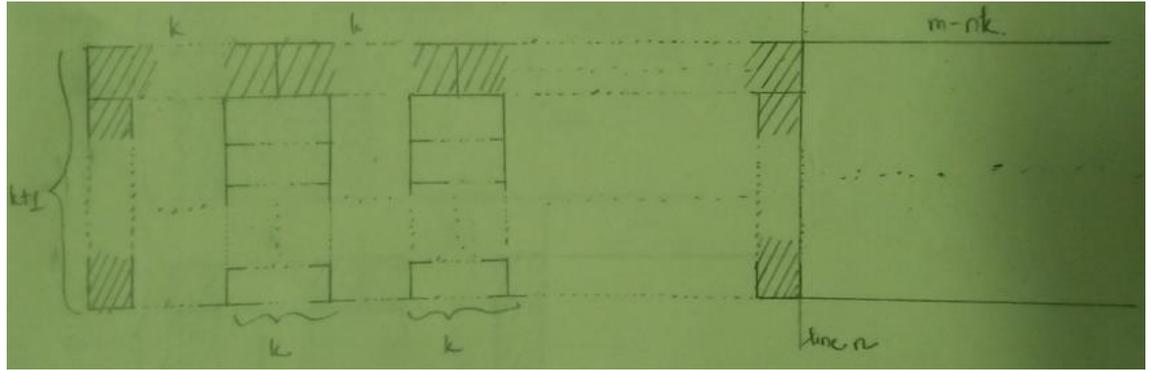

For the remaining untiled portion to the left of line $n$, since each line from 1 to $n-1$ is intersected by a different $k \times k$ square, hence the whole portion consists of $n-1$ $k \times k$ squares and $k-2$ vertical $k \times 1$ tiles.

Now we need to find the number of ways to arrange $n-1$ $k \times k$ squares and $k-2$ vertical tiles. Since all the arrangements will fill up the remaining space, we just need to pay attention to the number of ways to arrange the $k-2$ vertical tiles in the $n$ gaps created by the $n-1$ squares (some gaps may have no tiles while some gaps may have more than 1).

Let $x_1$ be the number of vertical tiles in the $1^{st}$ gap, $x_2$ be the number of vertical tiles in the $2^{nd}$ gap, …, $x_n$ be the number of vertical tiles in the $n^{th}$ gap.

*Then $x_1, x_2, \ldots, x_n$ are non-negative integers*

*and $x_1 + x_2 + x_3 + \ldots + x_n = k - 2$.*

*Each solution of $(x_1, x_2, \ldots, x_n)$ corresponds to one possible tiling.*

*Hence here we need to find the number of solutions $(x_1, x_2, \ldots, x_n)$*

*that satisfies the given conditions.*

*Each solution corresponds to a series of 0s and 1s :*

*00..00100…001…….00100…000, where*

*There are $n-1$ number 1 that acts as $n-1$ "partition walls" between the $n$ series of 0,*

*where the $1^{st}$ series has $x_1$ number 0s, the $2^{nd}$ series has $x_2$ number 0s, …*

*and the $n^{th}$ series has $x_n$ number 0s.*

*For example, a solution $(3,1,2)$ can be represented as 00010100.*

*There are a total of $k-2$ number 0s and $n-1$ number 1s in the series,*

*hence there are $k+n-3$ numbers in the series.*

*The number of distinct solutions corresponds to the number of arrangements of*

*0s and 1s in the series, which is also the number of ways to place $n-1$ number 1s*

*in the $k+n-3$ numbers in the series.*

*Thus the number of possible tilings (the number of solutions) is $^{k+n-3}C_{n-1}$.*

For the untiled portion to the right of line $n$, it is a $(k+1) \times (m-kn)$ board, and there are $a(k, m-kn)$ ways to tile.

Thus, in total there are $^{k+n-3}C_{n-1}$ **x a($k$,$m$-$kn$)** ways to tile the board.

- To account for all the cases, then we must consider all the cases of $n$ from 1 to $p-1$.
- The last case: when every line is intersected, then the top row is filled with $p$ horizontal $k \times 1$ tiles, and the rightmost tile is vertical, as follows.

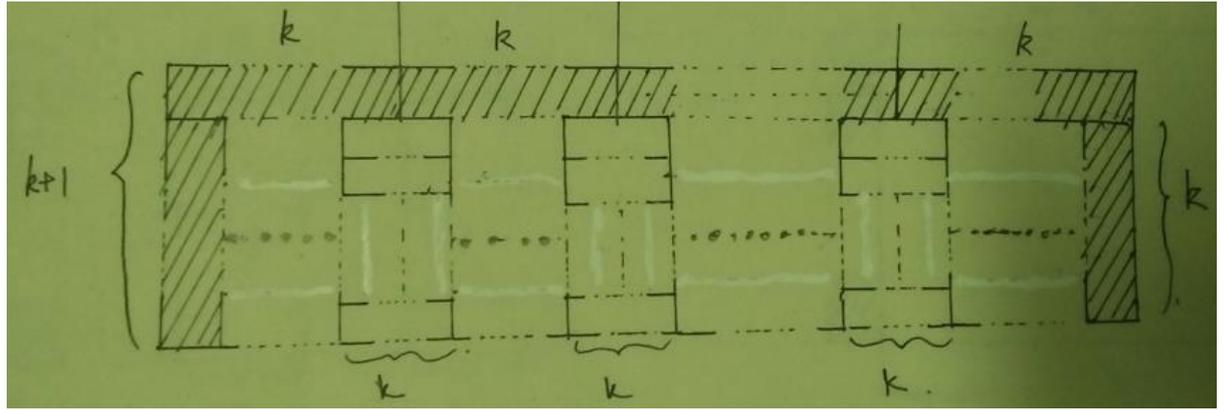

The remaining untiled portion consists of (p-1) k x k squares and k-2 vertical kx1 tiles. Using the above solution, substituting n for p, there are $^{k+p-3}C_{p-1}$ ways to arrange the k-2 vertical tiles in the p gaps created by the (p-1) squares, thus there are $^{k+p-3}C_{p-1}$ ways to tile the board.

In summary, for case **b**, adding up all the subcases for the different values of n, there are
$$a_{k,m-k} + {}^{k-1}C_1 a_{k,m-2k} + {}^{k}C_2 a_{k,m-3k} + {}^{k+1}C_3 a_{k,m-4k} + ... + {}^{k+p-5}C_{p-3} a_{k,2k} + {}^{k+p-4}C_{p-2} a_{k,k} + {}^{k+p-3}C_{p-1}$$
ways to tile the board.

e) When the leftmost column is covered with 1 horizontal and 1 vertical kx1 tile forming an "L" shape, this case is symmetrical to case **b**, hence there are also
$$a_{k,m-k} + {}^{k-1}C_1 a_{k,m-2k} + {}^{k}C_2 a_{k,m-3k} + {}^{k+1}C_3 a_{k,m-4k} + ... + {}^{k+p-5}C_{p-3} a_{k,2k} + {}^{k+p-4}C_{p-2} a_{k,k} + {}^{k+p-3}C_{p-1}$$
ways to tile the board.

Thus the number of ways to tile a (k+1) x m rectangular board using kx1 tiles, with m divisible by k, is

a(k,m)=a(k,pk)=a(k,m-k)+2(
$$a_{k,m-k} + {}^{k-1}C_1 a_{k,m-2k} + {}^{k}C_2 a_{k,m-3k} + {}^{k+1}C_3 a_{k,m-4k} + ... + {}^{k+p-5}C_{p-3} a_{k,2k} + {}^{k+p-4}C_{p-2} a_{k,k} + {}^{k+p-3}C_{p-1})$$

$\Rightarrow a_{k,m} = 3 a_{k,m-k} + 2({}^{k-1}C_1 a_{k,m-2k} + {}^{k}C_2 a_{k,m-3k} + {}^{k+1}C_3 a_{k,m-4k} + ... + {}^{k+p-5}C_{p-3} a_{k,2k} + {}^{k+p-4}C_{p-2} a_{k,k} + {}^{k+p-3}C_{p-1})$

### 4. Table of values

Using the formulae above, we can construct the table of values for a(k,m) with m divisible by k, or a(k,pk), the number of ways to tile a (k+1) x m rectangular board (k+1 rows, m columns) using kx1 tiles, for relatively small values of p and k. In this table, the first column represents values of k, while the first row represent values of p.

|   | 1 | 2 | 3 | 4 | 5 | 6 | 7 | 8 | 9 |
|---|---|---|---|---|---|---|---|---|---|
| 2 | 3 | 11 | 41 | 153 | 571 | 2131 | 7953 | 29681 | 110771 |
| 3 | 3 | 13 | 57 | 249 | 1087 | 4745 | 20713 | 90417 | 394691 |
| 4 | 3 | 15 | 75 | 371 | 1833 | 9057 | 44753 | 221137 | 1092699 |
| 5 | 3 | 17 | 95 | 521 | 2853 | 15629 | 85625 | 469105 | 2570031 |
| 6 | 3 | 19 | 117 | 701 | 4193 | 25101 | 150289 | 899825 | 5387479 |

## IV. Domino tiling of 4x*m* rectangular board
### 1. Domino tiling of 4x*m* rectangular board

Having tried and explored the cases of domino tilings of 2x*m* and 3x*m* rectangular boards, we will now move on to the case for a 4x*m* rectangular board.

For easy reference, label the columns, from left to right, 1,2,3,…, and the 4 rows, from top to bottom, A, B, C, D. The small square on column 3, row D, is referred to as D3, for instance.

Generally, we start from the bottom left corner of the table and consider 2 cases: when the orientation of the tile covering D1 (the "D1 tile") is either horizontal or vertical. Next we consider the same process for the adjacent small squares, going from left to right, "forcing" tiles whenever possible, in order to ensure that all cases are covered. This also ensures that no two cases overlap each other.

**Solution**

**For *m*=1,** it is obvious that there is only 1 way to tile a 4x1 rectangular board using 2 2x1 tiles. Thus **a(1)=1**.

**For *m*=2,** we have a 4x2 board. Consider 2 cases:

a) When the D1 tile is vertical:
   - If the B1 tile is also vertical, then for the remaining 4x1 board there is only 1 way to tile (Fig 1.1)
   - If the B1 tile is horizontal, then the rest of the tiles are "forced", and there is 1 way to tile (Fig 1.2).
b) When the D1 tile is horizontal: consider the C1 tile:
   - If the C1 tile is vertical, then the rest of the tiles are "forced", and there is 1 way (Fig 1.3).
   - If the C1 tile is horizontal, then for the remaining 2x2 board, there are 2 ways to tile (Fig 1.4).

In total, there are 1+1+1+2=5 ways to tile a 4x2 rectangular board using 4 2x1 tiles.

Thus **a(2)=5**.

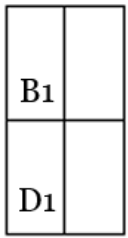 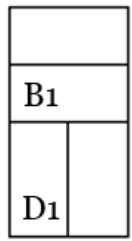 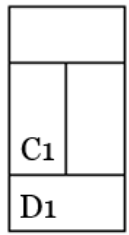 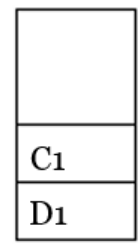

Fig 1.1  Fig 1.2  Fig 1.3  Fig 1.4

**For *m*=3,** we have a 4x3 board. Using the same process, considering the orientation of the D1, B1 and C1 tiles:

a) When the D1 and B1 tiles are both vertical, we are left with a 4x2 board (Fig 2.1). There are 5 ways to tile this board, as proven above.
b) When the D1 and C1 tiles are both horizontal, consider the B1 tile:
   - If the B1 tile is horizontal, then the A1 tile is "forced" as in Fig 2.2, and we are left with a 4x1 board. There is only 1 way to tile this board, as proven above.
   - If the B1 tile is vertical (Fig 2.3), consider the B2 tile:
     • If the B2 tile is horizontal, then the rest of the tiles are "forced" (Fig 2.4), and there is 1 way to tile the board.
     • If the B2 tile is vertical (Fig 2,5), we are left with a 4x1 board. There is only 1 way to tile this board.

Thus when the B1 tile is vertical, there are 2 ways to tile the board.

c) When the D1 tile is vertical and the B1 tile is horizontal, then the A1 tile is "forced". The board is left as in Fig 2.6. By rule of symmetry, the number of ways to tile is equal to the case in Fig 2.3, and there are 2 ways to tile the board.
d) When the D1 tile is horizontal and the C1 tile is vertical, the remaining tiles are "forced" as in Fig 2.7. There is only 1 way to tile the board.

Thus **a(3)=5+1+2+2+1=11.**

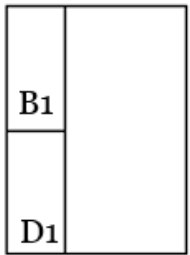 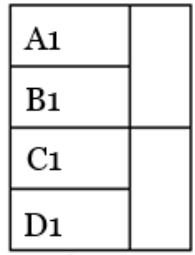 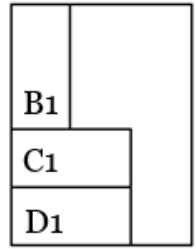 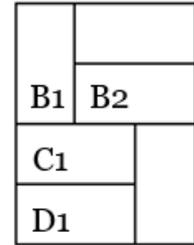
Fig 2.1  Fig 2.2  Fig 2.3  Fig 2.4

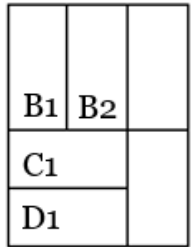 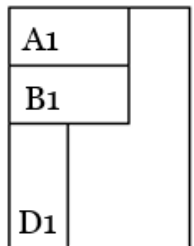 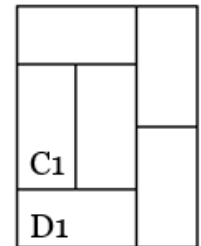
Fig 2.5  Fig 2.6  Fig 2.7

**For *m*=4,** we have a 4x4 board. Applying the same approach, we consider the 4 cases:

a) When the D1 and B1 tiles are both vertical, we are left with a 4x3 board, similar to Fig 2.1. There are 11 ways to tile this board, as proven above.
b) When the D1 and C1 tiles are both horizontal, consider the B1 tile:
   - If the B1 tile is horizontal, then the A1 tile is "forced", similar to Fig 2.2, and we are left with a 4x2 board. There are 5 ways to tile this board, as proven above.
   - If the B1 tile is vertical, similar to Fig 2.3, consider the B2 tile:
     • If the B2 tile is vertical, we are left with a 4x2 board. There are 5 ways to tile this board, as proven above.
     • If the B2 tile is horizontal, then the A2 tile is "forced" (Fig 3.1). The remaining board is exactly the same as in Fig 2.6, thus there are 2 ways to tile the board.

   Thus when the B1 tile is vertical, there are 5+2=7 ways to tile the board.

c) When the D1 tile is vertical and the B1 tile is horizontal, then the A1 tile is "forced", similar to Fig 2.6. As same as for *m*=2, by rule of symmetry, the number of ways to tile is equal to the case similar to Fig 2.3, when the B1 tile is vertical and there are 7 ways to tile the board.
d) When the D1 tile is horizontal and the C1 tile is vertical, the tiles are "forced" as in Fig 3.2. Consider the B2 tile:
   - If the B2 tile is vertical (Fig 3.3), we are left with a 4x2 board. There are 5 ways to tile this board, as proven above.

- If the B2 tile is horizontal, then the remaining tiles are "forced" as in Fig 3.4. There is only 1 way to tile the board.

Thus **a(4)=11+5+7+7+5+1=36**.

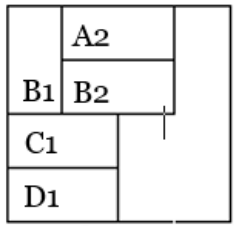
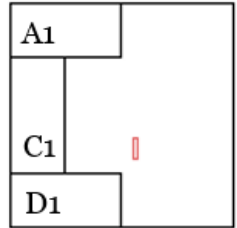
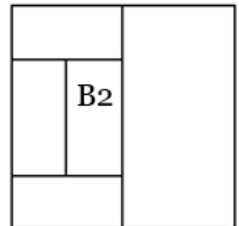
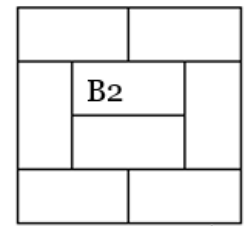

**For *m*≥5**, we apply the same approach and consider the 4 cases:

a) When the D1 and B1 tiles are both vertical, we are left with a 4x(*m*-1) board, similar to Fig 2.1. **There are a(*m*-1) ways to tile this board**.
b) When the D1 and C1 tiles are both horizontal, consider the B1 tile:
    - If the B1 tile is horizontal, then the A1 tile is "forced", similar to Fig 2.2, and we are left with a 4x(*m*-2) board. **There are a(*m*-2) ways to tile this board.**
    - If the B1 tile is vertical, similar to Fig 2.3, consider the B2 tile:
        • If the B2 tile is vertical, we are left with a 4x(*m*-2) board. There are a(*m*-2) ways to tile this board.
        • If the B2 tile is horizontal, then the A2 tile is "forced" (Fig 4.1). Consider the D3 tile:
            ➢ If the D3 tile is vertical, then we are left with a 4x(*m*-3) board. There are a(*m*-3) ways to tile this board.
            ➢ If the D3 tile is horizontal, then the C3 tile is "forced" (Fig 4.2). Continue to consider the B4 tile:
                ✓ If the B4 tile is vertical, we are left with a 4x(*m*-4) board. There are a(*m*-4) ways to tile this board.
                ✓ If the B4 tile is horizontal, then the A4 tile is "forced" (Fig 4.3). We will repeat this process in a "zig zag" manner until we reach the end of the 4x*m* board, when the last tile is "forced" (Fig 4.4)
    **Thus when the B1 tile is vertical, there are a(*m*-2)+a(*m*-3)+a(*m*-4)+…+a(1)+1 ways to tile the board.** This result stays true whether *m* is odd or even.

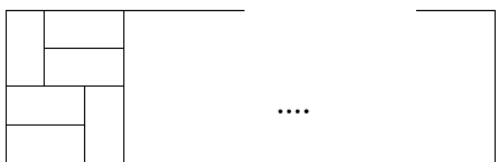
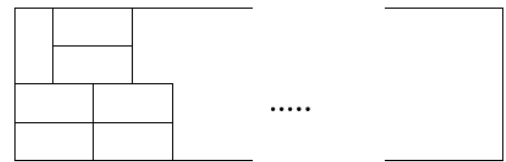
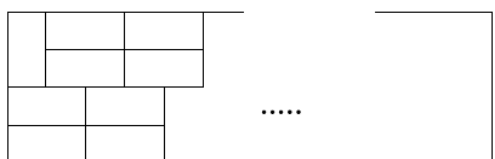
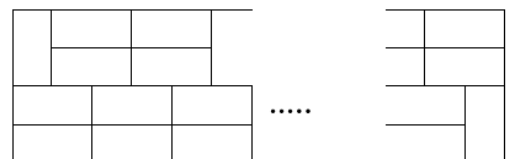

c) When the D1 tile is vertical and the B1 tile is horizontal, then the A1 tile is "forced", similar to Fig 2.6. By rule of symmetry, this case is similar to the case where D1 and C1 tiles are horizontal and the B1 tile is vertical. **Thus the number of ways in this case is also the same, which equals to a(*m*-2)+a(*m*-3)+a(*m*-4)+…+a(1)+1 ways.**

d) When the D1 tile is horizontal and the C1 tile is vertical, the tiles are "forced", similar to Fig 3.2. Consider the B2 tile:
   - If the B2 tile is vertical (similar to Fig 3.3), we are left with a 4x($m$-2) board. There are a($m$-2) ways to tile this board.
   - If the B2 tile is horizontal, then the tiles are "forced" as in Fig 4.5. Continue to consider tile B4:
     - If the B4 tile is vertical, then we are left with a 4x($m$-4) board. There are a($m$-4) ways to tile this board.
     - If the B4 tile is horizontal, then the tiles are "forced" as in Fig 4.6. We continue to repeat this process for tile B6 and beyond, "jumping" 2 columns at a time, until we reach the end of the board, when we cannot "jump" anymore. Here the results are slightly different depending on whether $m$ is odd or even:
       - If $m$ is odd, then by the time we reach the end of the board it would be like Fig 4.7. There is only 1 possible way to tile the remainder of the board.
         Thus for odd $m$, the number of ways to tile the board when the D1 tile is horizontal and the C1 tile is vertical is equal to **a($m$-2)+a($m$-4)+…+a(3)+1**.
       - If $m$ is even, then by the time we reach the end of the board it would be like Fig. 4.8. The rest of the tiles are "forced", contributing to 1 more way to tile the board. Thus for even $m$, the number of ways to tile the board when the D1 tile is horizontal and the C1 tile is vertical is equal to **a($m$-2)+a($m$-4)+…+a(2)+1**.

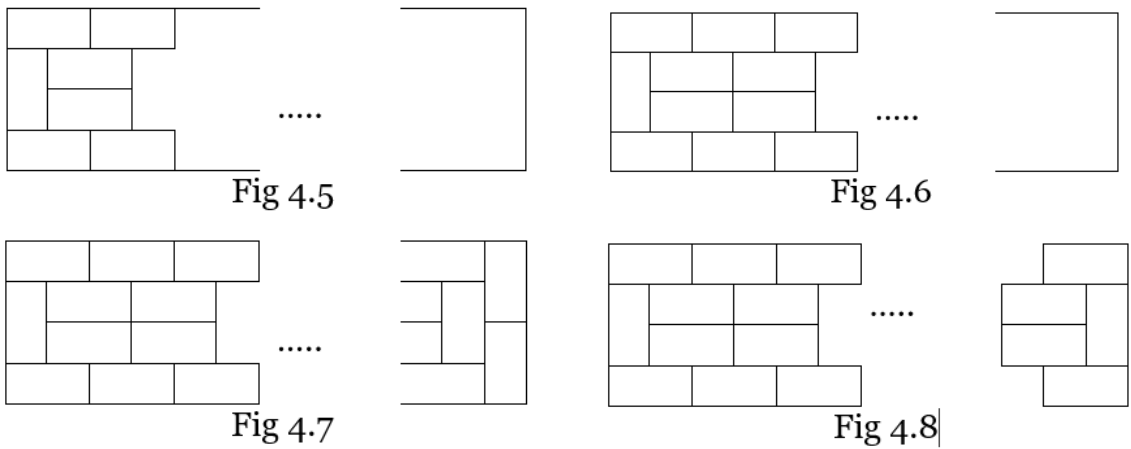

Fig 4.5   Fig 4.6

Fig 4.7   Fig 4.8

To summarize:

a) For even $m$, the total number of ways to tile the board is
   a($m$)= a($m$-1) + a($m$-2) + (a($m$-2) + a($m$-3) + a($m$-4) + … + a(2) + a(1) + 1)
                    + (a($m$-2) + a($m$-3) + a($m$-4) + … + a(2) + a(1) + 1)
                    + (a($m$-2) +           + a($m$-4) + … + a(2) +           1)
   Hence **a($m$)=a($m$-1)+4a($m$-2)+2a($m$-3)+3a($m$-4)+…+3a(2)+2a(1)+3**

b) For odd $m$, the total number of ways to tile the board is
   a($m$)= a($m$-1) + a($m$-2) + (a($m$-2) + a($m$-3) + a($m$-4) + … + a(3) + a(2) + a(1) + 1)
                     + (a($m$-2) + a($m$-3) + a($m$-4) + … + a(3) + a(2) + a(1) + 1)
                     + (a($m$-2) +           + a($m$-4) + … + a(3) +                1)
   Hence **a($m$)=a($m$-1)+4a($m$-2)+2a($m$-3)+3a($m$-4)+…+3a(3)+2a(2)+2a(1)+3**

Thus, when $m$ is odd:

a($m$) = a($m$-1) + 4a($m$-2) + 2a($m$-3) + 3a($m$-4) + 2a($m$-5) + … + 3a(3) + 2a(2) + 2a(1) + 3

            a($m$-2) = a($m$-3) + 4a($m$-4) + 2a($m$-5) + … +3a(3) + 2a(2) + 2a(1) + 3

⇨ a($m$) – a($m$-2) = a($m$-1) + 4a($m$-2) + a($m$-3) – a($m$-4)

⇨   a(*m*)         = a(*m*-1) + 5a(*m*-2) + a(*m*-3) – a(*m*-4)

When *m* is even, the same result is achieved using the same methods.

Thus the number of ways to tile an 4x*m* rectangular board using 2x1 tiles, a(*m*), is given by the formula

**a(1)=1, a(2)=5, a(3)=11, a(4)=36**

**a(*m*)=a(*m*-1)+5a(*m*-2)+a(*m*-3)-a(*m*-4)**

We will now attempt to construct a closed-form solution (formula of a(*m*) dependent only on *m*) from the recursive formula.

Using the recursive formula,

- When *m*=4, a(4) = a(3) + 5a(2) + a(1) – a(0)
  36 = 11 + 5x5 + 1 - a(0)
  a(0) = 1
- When *m*=3, a(3) = a(2) +5a(1) + a(0) – a(-1)
  11 = 5 + 5x1 + 1 - a(-1)
  a(-1) = 0

The characteristic equation is:

$$x^4 - x^3 - 5x^2 - x + 1 = 0$$

$$\Leftrightarrow x^2 - x - 5 - \frac{1}{x} + \frac{1}{x^2} = 0$$

$$\Leftrightarrow (x^2 + \frac{1}{x^2} + 2) - (x + \frac{1}{x}) - 7 = 0$$

$$\Leftrightarrow (x + \frac{1}{x})^2 - (x + \frac{1}{x}) - 7 = 0$$

Let $y = x + \frac{1}{x}$, then the equation becomes:

$$y^2 - y - 7 = 0$$

$$y = \frac{1 \pm \sqrt{29}}{2}$$

$$\Leftrightarrow x + \frac{1}{x} = \frac{1 \pm \sqrt{29}}{2}$$

$$\Leftrightarrow x^2 - \frac{1 \pm \sqrt{29}}{2} x + 1 = 0$$

$$\sqrt{\Delta} = \sqrt{\left(\frac{1 \pm \sqrt{29}}{2}\right)^2 - 4}$$

$$= \sqrt{\frac{14 \pm 2\sqrt{29}}{4}}$$

$$= \sqrt{\frac{(\sqrt{7+2\sqrt{5}} \pm \sqrt{7-2\sqrt{5}})^2}{4}}$$

$$= \frac{\sqrt{7+2\sqrt{5}} \pm \sqrt{7-2\sqrt{5}}}{2}$$

Let $A = \sqrt{7+2\sqrt{5}} + \sqrt{7-2\sqrt{5}}$

and $B = \sqrt{7+2\sqrt{5}} - \sqrt{7-2\sqrt{5}}$

For $x^2 - \frac{1+\sqrt{29}}{2}x + 1 = 0$,

$\sqrt{\Delta} = \frac{A}{2}$,

$$x_1 = \frac{1+\sqrt{29}+A}{4}$$

$$x_2 = \frac{1+\sqrt{29}-A}{4}$$

Note that the product of the two roots $x_1 x_2 = 1$

For $x^2 - \frac{1-\sqrt{29}}{2}x + 1 = 0$,

$\sqrt{\Delta} = \frac{B}{2}$,

$$x_3 = \frac{1-\sqrt{29}+B}{4}$$

$$x_4 = \frac{1-\sqrt{29}-B}{4}$$

Note that the product of the two roots $x_3 x_4 = 1$

Then the non-recurrence formula of a(m) has the form

$a(m) = ax_1^m + bx_2^m + cx_3^m + dx_4^m$

Sub $m = 0$: $a(0) = a + b + c + d$

$\Rightarrow a + b + c + d = 1$  (1)

Sub $m = 1$: $a(1) = ax_1 + bx_2 + cx_3 + dx_4$

$\Rightarrow a(1 + \sqrt{29} + A) + b(1 + \sqrt{29} - A) + c(1 - \sqrt{29} + B) + d(1 - \sqrt{29} - B) = 4$

$\Leftrightarrow (a + b - c - d)\sqrt{29} + A(a - b) + B(c - d) = 3$  (2)

Sub $m = -1$: $a(-1) = ax_1^{-1} + bx_2^{-1} + cx_3^{-1} + dx_4^{-1}$

$\Rightarrow ax_2 + bx_1 + cx_4 + dx_3 = 0$

$\Rightarrow a(1 + \sqrt{29} - A) + b(1 + \sqrt{29} + A) + c(1 - \sqrt{29} - B) + d(1 - \sqrt{29} + B) = 0$

$\Rightarrow (a + b - c - d)\sqrt{29} - A(a - b) - B(c - d) = -1$  (3)

From (2) and (3):

$(a + b - c - d)\sqrt{29} = 1$

$\Leftrightarrow (a + b) - (c + d) = \dfrac{\sqrt{29}}{29}$  (4)

From (1) and (4):

$a + b = \dfrac{29 + \sqrt{29}}{58}$  (5)

and $c + d = \dfrac{29 - \sqrt{29}}{58}$  (6)

Also, from (2) and (3):

$A(a - b) + B(c - d) = 2$  (7)

Sub $m = 2$: $a(2) = ax_1^2 + bx_2^2 + cx_3^2 + dx_4^2$

$\Rightarrow a(1 + \sqrt{29} + A)^2 + b(1 + \sqrt{29} - A)^2 + c(1 - \sqrt{29} + B)^2 + d(1 - \sqrt{29} - B)^2 = 80$

$\Leftrightarrow a(44 + 4\sqrt{29} + 2A + 2A\sqrt{29}) + b(44 + 4\sqrt{29} - 2A - 2A\sqrt{29})$

$+ c(44 - 4\sqrt{29} + 2B - 2B\sqrt{29}) + d(44 - 4\sqrt{29} - 2B + 2B\sqrt{29}) = 80$

$\Leftrightarrow 44(a + b + c + d) + 4\sqrt{29}(a + b - c - d)$

$+ 2A(a - b) + 2B(c - d) + 2A\sqrt{29}(a - b) - 2B\sqrt{29}(c - d) = 80$

$\Leftrightarrow 44 + 4 + 4 + 2A\sqrt{29}(a - b) - 2B\sqrt{29}(c - d) = 80$

$\Leftrightarrow A\sqrt{29}(a - b) - B\sqrt{29}(c - d) = 14$

$\Leftrightarrow A(a - b) - B(c - d) = \dfrac{14\sqrt{29}}{29}$  (8)

From (7) and (8):

$A(a - b) = \dfrac{29 + 7\sqrt{29}}{29}$  (9)

and $B(c - d) = \dfrac{29 - 7\sqrt{29}}{29}$  (10)

From (5) and (9):

$$a = (\frac{29+\sqrt{29}}{116}) + \frac{29+7\sqrt{29}}{58A} = \frac{\sqrt{29}(1+\sqrt{29})}{116} + \frac{1}{A}(\frac{29+7\sqrt{29}}{58})$$

$$b = (\frac{29+\sqrt{29}}{116}) - \frac{29+7\sqrt{29}}{58A} = \frac{\sqrt{29}(1+\sqrt{29})}{116} - \frac{1}{A}(\frac{29+7\sqrt{29}}{58})$$

*From (6) and (10):*

$$c = (\frac{29-\sqrt{29}}{116}) + \frac{29-7\sqrt{29}}{58B} = \frac{-\sqrt{29}(1-\sqrt{29})}{116} + \frac{1}{B}(\frac{29-7\sqrt{29}}{58})$$

$$b = (\frac{29-\sqrt{29}}{116}) - \frac{29+7\sqrt{29}}{58B} = \frac{-\sqrt{29}(1-\sqrt{29})}{116} - \frac{1}{B}(\frac{29-7\sqrt{29}}{58})$$

For simplification, let

$C = 1 + \sqrt{29}$

$D = 1 - \sqrt{29}$

$E = \dfrac{29 + 7\sqrt{29}}{58}$

$F = \dfrac{29 - 7\sqrt{29}}{58}$

$G = \dfrac{\sqrt{29}}{116}$

Thus the non-recursive formula for a(m) is

$$a(m) = (CG + \frac{E}{A})(\frac{A+C}{4})^m + (CG - \frac{E}{A})(\frac{C-A}{4})^m + (-DG + \frac{F}{B})(\frac{B+D}{4})^m + (-DG - \frac{F}{B})(\frac{D-B}{4})^m$$

As we reach this conclusion, we ask ourselves the same question that we did above: is it possible to replace 2 by 3, by *k* while maintaining the configuration of the board relative to the tile? In other words, is it possible to find the number of tilings of a 6x*m* rectangular board with 3x1 trominoes? A *(2k) x m* rectangular board with *k*x1 tiles?

Unfortunately, we have attempted it but ended up with little success.

## 2. Table of values

Using the above formula, we are able to construct a table of values for a(*m*), the number of ways to tile a 4x*m* rectangular board using 2x1 tiles, with *m* relatively small.

| m | 1 | 2 | 3 | 4 | 5 | 6 | 7 | 8 | 9 | 10 | 11 | 12 | 13 | 14 |
|---|---|---|---|---|---|---|---|---|---|---|---|---|---|---|
| a(m) | 1 | 5 | 11 | 36 | 95 | 281 | 781 | 2245 | 6336 | 18061 | 51205 | 145601 | 413351 | 1174500 |

# V. Tatami tiling
## 1. Definition

A tatami tiling is a domino tiling which satisfies the "tatami property" – that no four tiles meet at a point.

Examples of valid tatami tilings are as follows:

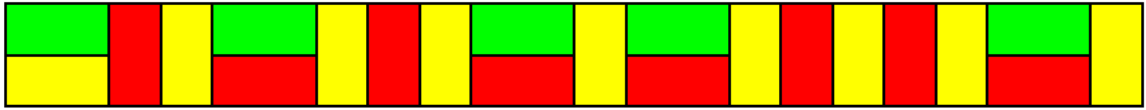

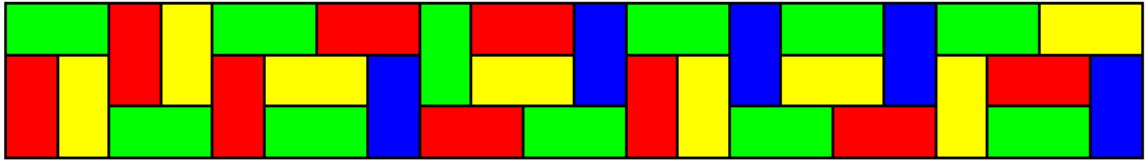

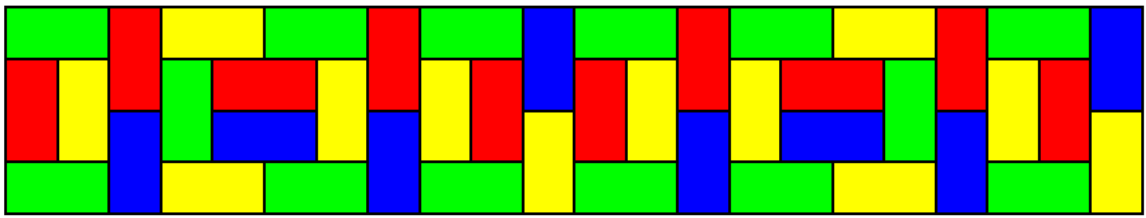

It is important to note that in this report, we are only concerned with the number of unrestricted tatami tilings, in which tilings are counted irrespective of internal symmetry. Tilings that match each other after rotations and/or reflections are counted with their multiplicity (rotations and reflections are considered distinct tilings).

## 2. Tatami tiling of 2x*m* rectangular board

Generally speaking, the approach to the problem of tatami tiling is same as the ordinary tiling problem. However, unlike ordinary tiling, tatami tiling has certain restrictions that must be taken into account when counting.

We will start with the simplest case: the number of tatami tilings of a 2x*m* rectangular board, a(*m*).

**Solution**

It is obvious that **a(1)=1, a(2)=2**.

For *m*=3, there are 3 possible ways of tiling (as follows). Thus **a(3)=3**.

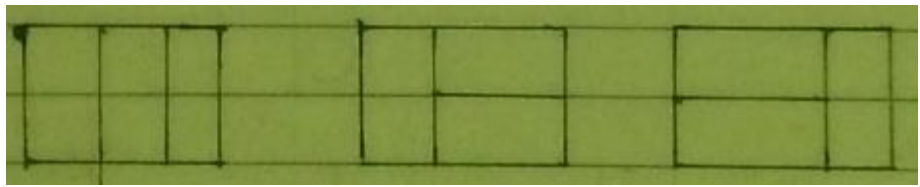

For *m*≥4, consider 2 cases involving the tiling of the leftmost column:

    a) When the leftmost column is covered with a vertical 2x1 domino:

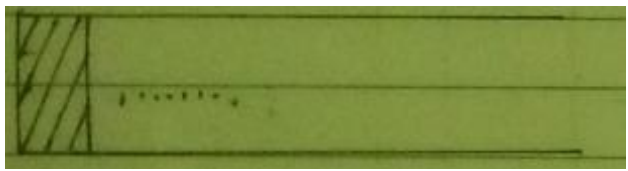

In this case, the first tile causes no restriction on the tiling of the remaining 2x1 dominoes (the next tile can either be a vertical or horizontal tile). Thus we are left with a 2x($m$-1) board, and there are a($m$-1) ways to tile.

b) When the leftmost column is covered with 2 horizontal 2x1 dominoes, then to preserve the tatami property, the next tile is "forced" into place as follows:

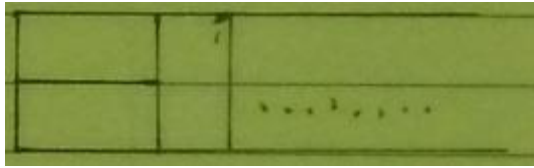

Since the tile on the 3<sup>rd</sup> column is vertical, the tiles on the 4<sup>th</sup> column can either be vertical or horizontal (there are no restrictions). Thus we are left with a 2x($m$-3) board, and there are a($m$-3) ways to tile.

In summary, the number of unrestricted tatami tilings of a 2x$m$ rectangular board is

**a($m$)=a($m$-1)+a($m$-3)** for $m \geq 4$

with a(1)=1, a(2)=2 and a(3)=3

## 3. Tatami tiling of 3x$m$ rectangular board

The case for the 2x$m$ rectangular board was fairly simple. How about for the 3x$m$ rectangular board? It is clearly more difficult, however, thanks to the restriction of the tatami property, we are able to reduce the huge number of cases into a few number of possible cases based on the tiling of the leftmost column.
**Solution**

We will move on to a slightly more difficult case: the number of tatami tilings of a 3x$m$ rectangular board, a($m$).

It is obvious that for any possible tiling to occur, 3x$m$ must be divisible by 2, hence $m$ must be divisible by 2. Thus **a($m$)=0 for odd $m$**.

We consider the case when $m$ is even.

Generally, there are 3 cases involving the tiling of the leftmost column, as follows.

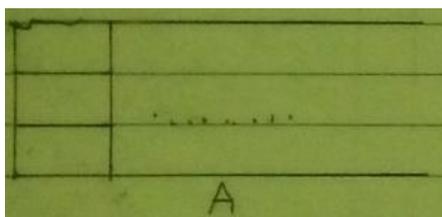 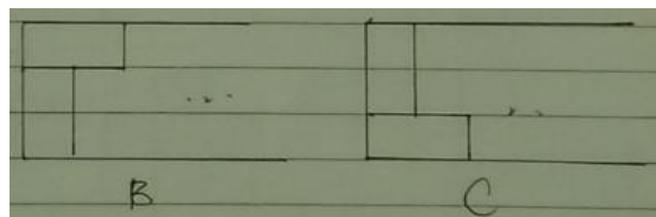

For case A, when $m$>2 there is no possible tiling of the third column that satisfies no 4 tiles meet at one corner. Thus for $m$>2 then any possible tiling must either belong to case B or case C.

Denote B($m$) and C($m$) as the number of tatami tilings of a 3x$m$ rectangular board with the leftmost column tiled like case B and C respectively.

Then for $m$>2, a($m$)=B($m$)+C($m$).

**For $m$=2**, we have a 3x2 board. Similar to the above case, there are 3 ways to tile. Thus **a(2)=3**.

**For $m$=4**, we have a 3x4 board. Let's consider 2 cases mentioned above:

a) When the leftmost column is tiled like case B, then we consider the tiling of the third tile:

- If the third tile is vertical, then to preserve the tatami property the other tiles are "forced" into place as follows, and there is only 1 way to tile.

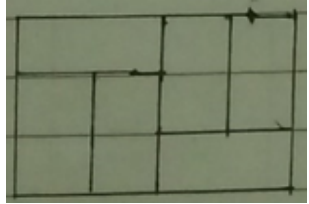

- If the third tile is horizontal, there is also only 1 way to tile to preserve the tatami property as follows:

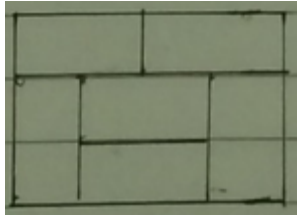

Thus for this case there are 2 ways to tile.
b) When the leftmost column is tiled like case C, it is symmetrical to case B, hence there are also 2 ways to tile.

In summary, the number of tatami tilings of a 3x4 board is **a(4)=2+2=4**.

**For *m*=6**, we have a 3x6 board. Similarly, consider 2 cases:

a) When the leftmost column is tiled like case B, then we consider the tiling of the third tile:
- If the third tile is vertical, then to preserve the tatami property the next tiles are "forced" into place as follows:

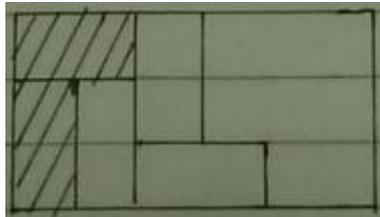

Here we consider the sixth tile, whether it's horizontal or vertical, and in both cases the other tiles are "forced" into place to preserve the tatami property, yielding 1 possible tiling each.

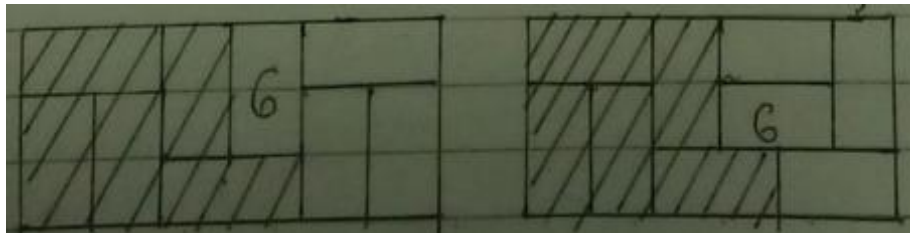

- If the third tile is horizontal, then the remaining tiles are "forced" place to preserve the tatami property. Thus there is 1 way to tile.

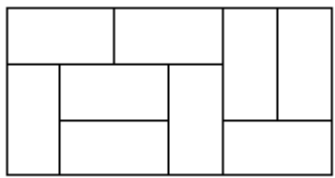

Thus for this case there are 3 ways to tile.

b) When the leftmost column is tiled like case C, it is symmetrical to case B, hence there are also 3 ways to tile.

In summary, the number of tatami tilings of a 3x6 board is **a(6)=3+3=6**.

**For *m*≥8**, using the same approach, consider the 2 cases:

a) When the leftmost column is tiled like case B, then we consider the tiling of the third tile:
   - If the third tile is vertical, then to preserve the tatami property the next tiles are "forced" into place as follows:

   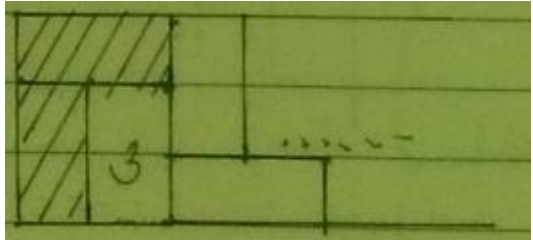

   The remaining untiled portion is exactly the same as a 3x(*m*-2) rectangular board with the leftmost column tiled like case C.
   Hence there are C(*m*-2) ways to tile.
   - If the third tile is horizontal, then due to the tatami property the next tiles are "forced":

   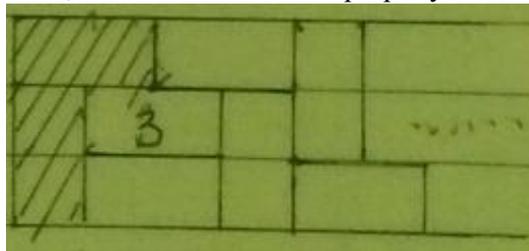

   The remaining untiled portion is exactly the same as a 3x(*m*-4) rectangular board with the leftmost column tiled like case C.
   Hence there are C(*m*-4) ways to tile.

   Thus for this case, there are **C(*m*-2)+C(*m*-4)** possible tilings.
b) When the leftmost column is tiled like case C, it is symmetrical to case B, hence in this case there are **B(*m*-2)+B(*m*-4)** possible tilings.

In summary, the number of tatami tilings of a 3x*m* board is

**a(*m*)**=C(*m*-2)+C(*m*-4)+B(*m*-2)+B(*m*-4)

=(B(*m*-2)+C(*m*-2))+(B(*m*-4)+C(*m*-4))

=**a(*m*-2)+a(*m*-4)** (since *m*≥8 so *m*-4>2)

Thus the formula for a(*m*) is

**a(*m*)=0 for odd *m***

**a(*m*)=a(*m*-2)+a(*m*-4) for *m* even and *m*≥8**

**with a(2)=3, a(4)=4 and a(6)=6**

## 4. Tatami tiling of 4x*m* rectangular board

As it turns out, the solution for the 4x*m* rectangular board is strikingly similar to the 3x*m* rectangular board, as we are also able to restrict the different ways of tiling of the leftmost column due to the tatami property.

**Solution**

Let's explore the case even further: is it possible to use the same approach to find the number of tatami tilings of a 4x*m* rectangular board, a(*m*)?

Generally, there are 5 cases involving the tiling of the leftmost column, as follows.

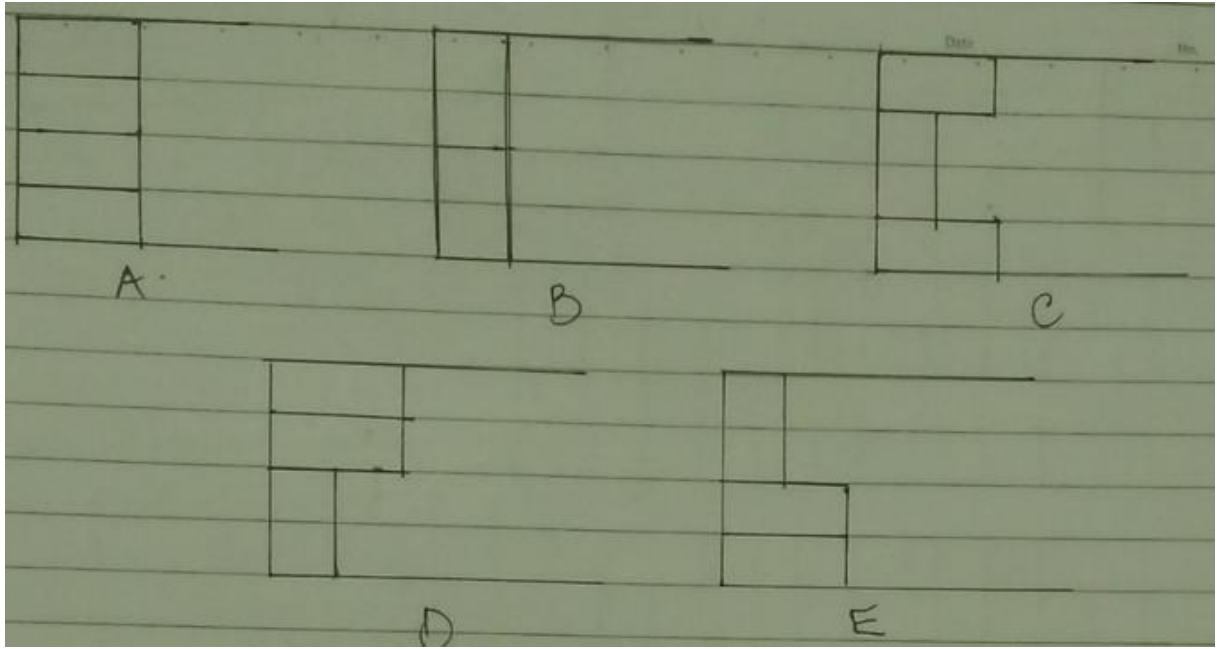

For case A, when *m*>2 there is no possible tiling of the third column that satisfies no 4 tiles meet at one corner.

Similarly for case D, there is no possible tiling of the 4th column that satisfies the tatami property (as follows). Case E is symmetrical to case D, thus it has the same property.

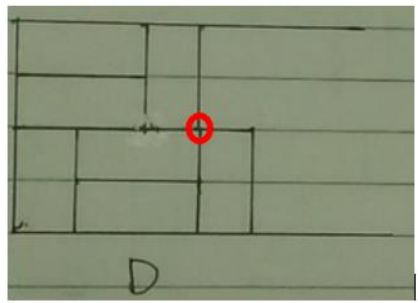

Thus for *m*>3, any possible tatami tiling of the 4x*m* rectangular board must either have the left most column tiled like case B or case C.

Denote B(*m*) and C(*m*) as the number of tatami tilings of a 4x*m* rectangular board with the leftmost column tiled like case B and C respectively.

Then for *m*>3, a(*m*)=B(*m*)+C(*m*).

However, this time, we will consider the general cases first, as the approach is useful to quickly solving the trivial cases.

**For *m* sufficiently large**, using the same approach, consider the 2 cases:

a) When the leftmost column is tiled like case B, then to preserve the tatami property the next tiles are "forced" into place as follows:

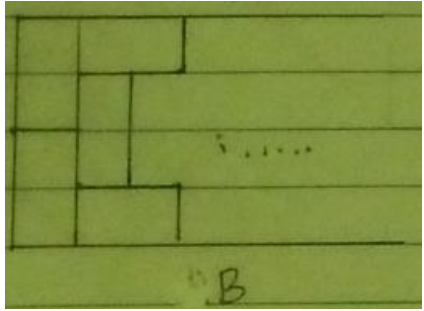

Now we consider 2 cases involving the orientation of the 6th tile:
- If the 6th tile is vertical, then to preserve the tatami property the next tiles are "forced" into place as follows (case B1).

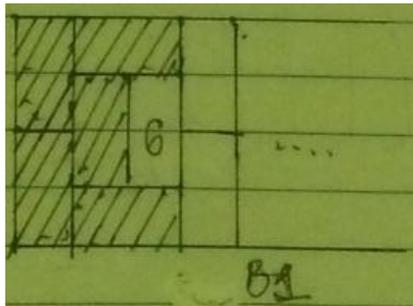

The remaining untiled portion is exactly the same as a 4x($m$-3) rectangular board with the leftmost column tiled like case B.

Hence there are B($m$-3) ways to tile.

- If the 6th tile is horizontal, then due to the tatami property the next tiles are "forced" (case B2):

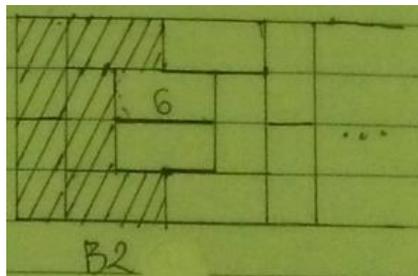

The remaining untiled portion is exactly the same as a 4x($m$-5) rectangular board with the leftmost column tiled like case B.

Hence there are B($m$-5) ways to tile.

Thus for this case, there are **B($m$-3)+B($m$-5)** possible tilings.

b) When the leftmost column is tiled like case C, then consider the orientation of the 4th tile:
- If the 4th tile is vertical, then to preserve the tatami property the next tiles are "forced" into place as follows (case C1).

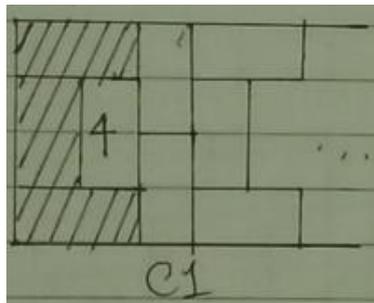

The remaining untiled portion is exactly the same as a 4x($m$-3) rectangular board with the leftmost column tiled like case C.

Hence there are C($m$-3) ways to tile.

- If the 4[th] tile is horizontal, then due to the tatami property the next tiles are "forced" (case C2):

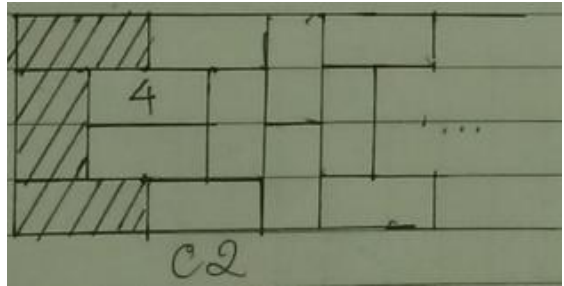

The remaining untiled portion is exactly the same as a 4x($m$-5) rectangular board with the leftmost column tiled like case C.

Hence there are C($m$-5) ways to tile.

Thus for this case, there are **C($m$-3)+C($m$-5)** possible tilings.

In summary, the number of tatami tilings of a 4x$m$ board is

**a($m$)=B($m$-3)+B($m$-5)+C($m$-3)+C($m$-5)**

=(B($m$-3)+C($m$-3))+(B($m$-5)+C($m$-5))

**=a($m$-3)+a($m$-5)**

This identity holds true for $m$-5>3, or $m$≥9.

Thus we will now find out the number of tatami tilings for $m$ from 1 to 8.

It is obvious that **a(1)=1**.

**For $m$=2 and $m$=3**, we have a 4x2 board and a 4x3 board and there are 4 possible tatami tilings for each case as shown:

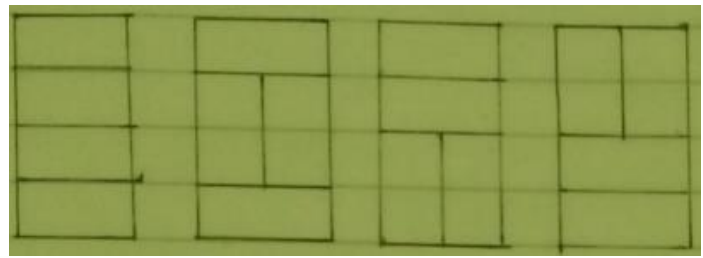

Tatami tiling cases for $m$=2

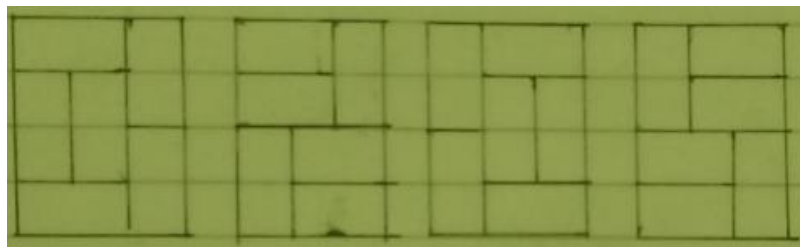

Tatami tiling cases for $m$=3

Thus **a(2)=a(3)=4**.

**For $m$=4,5,6**, the leftmost column can only be tiled like case B or case C. By considering each case and trying to "force" tiles into place to preserve the tatami property like in the general case, we can easily find out that for *m*=4 there are 2 possible tilings while for *m*=5 and 6 there are 3 possible tatami tilings for each case:

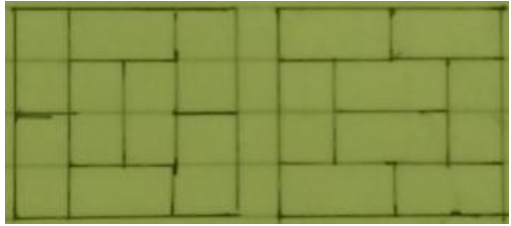

Tatami tiling cases for *m*=4

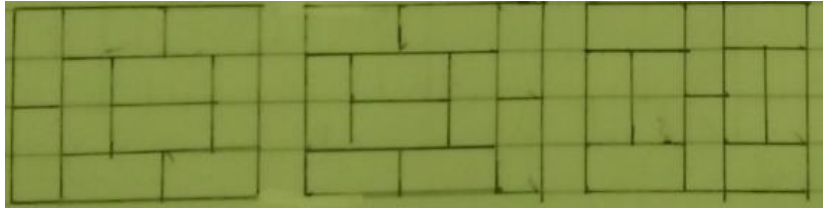

Tatami tiling cases for *m*=5

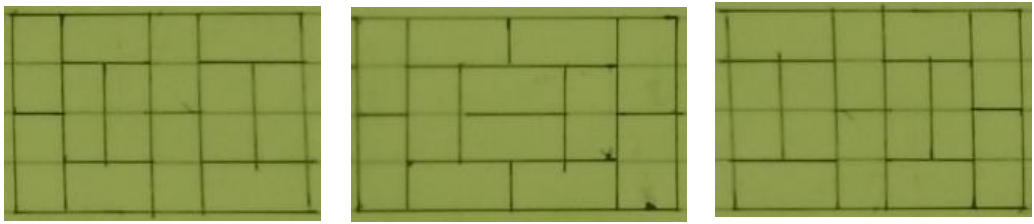

Tatami tiling cases for *m*=6

Thus **a(4)=2, a(5)=a(6)=3**.

**For *m*=7**, consider cases regarding the tiling of the leftmost column:

a) If the leftmost column is tiled like case B, then we are left with a 4x6 board. Out of the 3 possible tatami tilings for *m*=6, only one maintains the tatami property for *m*=7. Thus there is 1 way to tile.

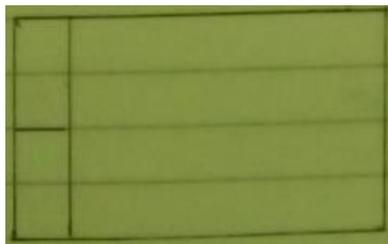

b) If the leftmost column is tiled like case C1, then we are left with a 4x5 board. Out of the 3 possible tatami tilings for *m*=5, only one maintains the tatami property for *m*=7. Thus there is 1 way to tile.

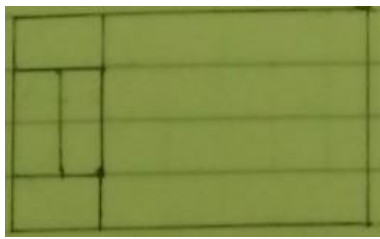

c) If the leftmost column is tiled like case C2, then we are left with a 4x3 board. Out of the 4 possible tatami tilings for *m*=3, only one maintains the tatami property for *m*=7. Thus there is 1 way to tile.

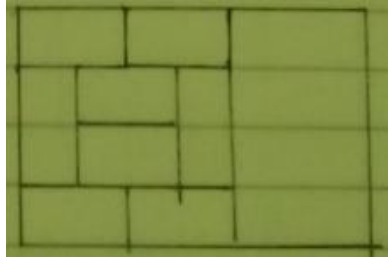

Thus the number of possible tatami tilings for a 4x7 rectangular board is **a(7)=3**.

**For *m*=8**, we apply the same approach as *m*=7, consider the tiling of the leftmost column:

a) If the leftmost column is tiled like case B, then we are left with a 4x7 board. Out of the 3 possible tatami tilings for *m*=7, only 2 maintains the tatami property for *m*=8. Thus there are 2 ways to tile.
b) If the leftmost column is tiled like case C1, then we are left with a 4x6 board. Out of the 3 possible tatami tilings for *m*=6, only 2 maintains the tatami property for *m*=8. Thus there are 2 ways to tile.

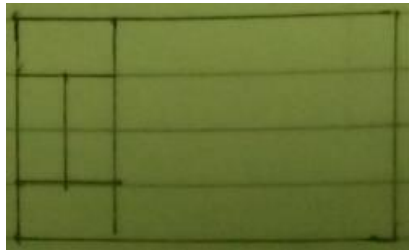

c) If the leftmost column is tiled like case C2, then we are left with a 4x4 board. Out of the 2 possible tatami tilings for *m*=4, only one maintains the tatami property for *m*=8. Thus there is 1 way to tile.

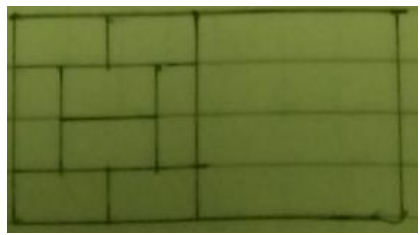

Thus the number of possible tatami tilings for a 4x8 rectangular board is **a(8)=2+2+1=5**.

We can conclude that the formula for a(*m*) is

$$a(m)=a(m-3)+a(m-5) \text{ for } m \geq 9$$

with **a(1)=1, a(2)=a(3)=4, a(4)=2, a(5)=a(6)=a(7)=3 and a(8)=5.**

## 5. Table of values

Based on the formulae deduced above, we can construct a table of values for a(*m,n*), the number of unrestricted tatami tilings of a *mxn* rectangular board, for *m*=2,3,4 and relatively small values of *n*. In this table, the first column represents values of *m*, while the first row represent values of *n*.

|   | 1 | 2 | 3 | 4 | 5 | 6 | 7 | 8 | 9 | 10 | 11 | 12 | 13 | 14 | 15 |
|---|---|---|---|---|---|---|---|---|---|----|----|----|----|----|----|
| 2 | 1 | 2 | 3 | 4 | 6 | 9 | 13 | 19 | 28 | 41 | 60 | 88 | 129 | 189 | 277 |
| 3 | 0 | 3 | 0 | 4 | 0 | 6 | 0 | 10 | 0 | 16 | 0 | 26 | 0 | 42 | 0 |
| 4 | 1 | 4 | 4 | 2 | 3 | 3 | 3 | 5 | 5 | 6 | 8 | 8 | 11 | 13 | 14 |

# Conclusion

Through this paper, we have managed to produce a simpler solution based on recursive sequence for 3 different configurations of the board. We have also solved the simple cases for tatami tiling. In conclusion, we hope to have shed light on this famous problem with a unique approach and made it more available and understandable to students. Our attempt uses mainly Combinatorics which most students can comprehend. However, it is also worth notice that our approach works on a case-by-case basis and we have failed to deliver a solution that works for the general case. Our approach can be used to solve many other problems in Combinatorics, especially counting problems, which we wish to further research in the future.

# Acknowledgements